\documentclass[aap,MSNbibl,seceqn,dvips]{arximspdf}

\doi{10.1214/12-AAP881} 
\volume{23}
\issue{4}
\pubyear{2013}
\firstpage{1584}
\lastpage{1628}

\makeatletter

\newcommand{\rrVert}{\Vert}
\newcommand{\rrvert}{\vert}
\newcommand{\llVert}{\Vert}
\newcommand{\llvert}{\vert}

\newcommand{\cal}{\mathcal}

\newcommand{\define}{\triangleq}
\newcommand{\ind}{\mathbb{I}}
\newcommand{\sF}{{\cal F}}
\newcommand{\sG}{{\cal G}}
\newcommand{\sH}{{\cal H}}
\newcommand{\frakC}{{\mathfrak C}}
\newcommand{\calE}{{\cal E}}
\newcommand{\sE}{{\cal E}}
\newcommand{\frakE}{{\mathfrak E}}

\newcommand{\R}{\mathbb{R}}
\renewcommand{\P}{\mathbb{P}}
\newcommand{\Q}{\mathbb{Q}}
\newcommand{\E}{\mathbb{E}}
\newcommand{\N}{\mathbb{N}}
\newcommand{\Nbar}{\overline{\N}}
\newcommand{\Xbar}{\overline{X}}
\newcommand{\Ybar}{\overline{Y}}
\newcommand{\Zbar}{\overline{Z}}
\newcommand{\sigmahat}{\widehat{\sigma}}
\newcommand{\What}{\widehat{W}}

\newcommand{\Omegahat}{\widehat{\Omega}}
\newcommand{\sFhat}{\widehat{\sF}}
\newcommand{\Phat}{\widehat{\P}}
\newcommand{\Xhat}{\widehat{X}}
\newcommand{\Yhat}{\widehat{Y}}

\newcommand{\Mhat}{\widehat{M}}
\newcommand{\Zhat}{\widehat{Z}}
\newcommand{\omegabar}{\overline{\omega}}
\newcommand{\Qhat}{\widehat{Q}}
\newcommand{\Chat}{\widehat{C}}
\newcommand{\Ahat}{\widehat{A}}
\newcommand{\Bhat}{\widehat{B}}
\newcommand{\ebar}{\overline{e}}
\newcommand{\sB}{{\cal B}}
\newcommand{\Gammahat}{\widehat{\Gamma}}
\newcommand{\sO}{\cal{O}}
\newcommand{\Bbar}{\overline{B}}
\newcommand{\Cbar}{\overline{C}}

\newtheorem{theorem}{Theorem}[section]
\newtheorem{corollary}[theorem]{Corollary}
\newtheorem{proposition}[theorem]{Proposition}
\newtheorem{lemma}[theorem]{Lemma}

\newproclaim{definition}[theorem]{Definition}
\newproclaim{example}[theorem]{Example}
\newproclaim{remark}[theorem]{Remark}

\makeatother

\begin{document}
\begin{frontmatter}

\title{Mimicking an It\^{o} process by
a solution of a~stochastic differential equation}
\runtitle{Mimicking an It\^{o} process}

\begin{aug}
\author[A]{\fnms{Gerard} \snm{Brunick}\thanksref{t1}\ead[label=e1]{brunick@pstat.ucsb.edu}}
\and
\author[B]{\fnms{Steven} \snm{Shreve}\corref{}\thanksref{t2}\ead[label=e2]{shreve@andrew.cmu.edu}}
\runauthor{G. Brunick and S. Shreve}
\affiliation{University of California and Carnegie Mellon University}
\address[A]{Department of Statistics\\
\quad and Applied Probability\\
University of California, Santa Barbara\\
Santa Barbara, California 93106\\
USA\\
\printead{e1}} 
\address[B]{Department of Mathematical Sciences\\
Carnegie Mellon University\\
Pittsburgh, Pennsylvania 15213\\
USA\\
\printead{e2}}
\end{aug}

\thankstext{t1}{Supported in part by the NSF Grant DMS-04-04682.}

\thankstext{t2}{Supported in part by the NSF Grants DMS-04-04682
and DMS-09-03475.}

\received{\smonth{10} \syear{2011}}
\revised{\smonth{5} \syear{2012}}

%
\begin{abstract}
Given a multi-dimensional It\^o process whose drift and diffusion terms
are adapted processes, we construct a weak solution to a stochastic
differential equation that matches the distribution of the It\^o
process at each fixed time. Moreover, we show how to match the
distributions at each fixed time of functionals of the It\^o process,
including the running maximum and running average of one of the
components of the process. A consequence of this result is that a wide
variety of exotic derivative securities have the same prices when the
underlying asset price is modeled by the original It\^o process or the
mimicking process that solves the stochastic differential equation.
\end{abstract}

%

\pdfkeywords{60G99, 60H10, 91G20, Ito process,
stochastic differential equation,
derivative security pricing,
stochastic volatility models}

\end{frontmatter}

\section{Introduction}\label{SIntro}

We construct a process that mimics certain
properties of a given It\^o process, but is
simpler in the sense that the mimicking process
solves a stochastic differential equation (SDE),
while the It\^o process may have drift and
diffusion terms that are themselves stochastic processes.
This work is motivated by the problem of model
calibration in finance. The financial engineer
would like to identify a class of models for
an underlying asset price that is
flexible enough to allow for calibration to a wide
range of possible market prices of derivative securities
on that asset. The result of this paper shows the extent to which
sophisticated models are no more powerful
for calibration purposes than an SDE
for the underlying asset price.

Our results are closely related to Krylov~\cite{Kry3} and
Gy\"ongy~\cite{G}. Krylov~\cite{Kry3} calls the measure that
records the average amount of time that
an It\^o process $X$ spends in each Borel set
before being killed at the first jump of an
independent Poisson
process with intensity $\lambda$ the
\textit{Green $\lambda$-measure} of $X$.
Given an It\^o process with\vadjust{\goodbreak}
bounded drift and bounded, uniformly
positive-definite covariance,
Krylov~\cite{Kry3} constructs a
process with the same Green $\lambda$-measure
which solves a
time-independent diffusion equation. Krylov
further asserts that it
is possible to construct a process that solves
a time-dependent
diffusion equation and matches the
one-dimensional marginal
distributions of such an It\^o process.
Gy\"ongy~\cite{G} provides a
proof of Krylov's assertion and shows that
the drift and
covariance in the diffusion
equation solved by the
mimicking process may be interpreted as the expected value of the
It\^o process's instantaneous drift and covariance conditioned on its
level.
See also Klebaner~\cite{Kleb} for a
related argument
based on semimartingale local time.

Gy\"ongy~\cite{G} was rediscovered by the mathematical finance
community in the context of local volatility models.
Dupire~\cite{Dup1} studies the European option prices generated by
a model in which the risk-neutral dynamics of the price process
satisfy a time-dependent diffusion equation
(see also Derman and Kani~\cite{DK1} for a discrete-time treatment
of this topic). These models are now
known as \textit{local volatility models}, and the diffusion coefficient
of the
log-price process is known as the \textit{local volatility surface}.
Dupire~\cite{Dup1} shows that it is possible to construct a local
volatility model that
is consistent with a given set of European option prices when that set
of prices is sufficiently smooth as a function of maturity and strike,
and he shows how the local volatility surface may be implied directly
from the call prices. Local volatility models have proven popular
with practitioners because they allow for calibration to a wide range
of European option prices. Dupire~\cite{Dup1} does not find the
dynamics of a local volatility model to be particularly plausible;
however, he asserts that ``the market prices European options
\textit{as if} the process was this diffusion.'' In effect, the local
volatility model mimics the European option prices of some more
complicated market process, and this is equivalent to matching the
one-dimensional marginal distributions of
that process under the equivalent martingale
probability measure (also call the
\textit{risk-neutral measure}) used for pricing.

In~\cite{Dup2}, Dupire extends~\cite{Dup1} to
study the local volatility surface that is implied
not by market prices of options but by prices
generated from a stochastic volatility model. Using
infinitesimal calendar and butterfly spreads, he
presents a financial argument that the square of
the local volatility function is the expected value
of the instantaneous squared stochastic volatility
conditioned on the level of the underlying asset
price, essentially recovering Gy\"ongy's result,
albeit in a nonrigorous fashion. Following this
development, the Gy\"ongy--Dupire formula has found
several applications in finance. For example,
Gatheral~\cite{Gat} uses it to compare the
properties of a number of stochastic volatility
models, and Antonov and Misirpashaev~\cite{AM}
and Piterbarg~\cite{Pit2,Pit4,Pit5} combine it with
parameter averaging techniques to produce pricing
approximations
based on approximations of the second conditional
expectation appearing in (\ref{2.18}) below,
a special case of (\ref{GD}) in our main result.

Brigo and Mecurio~\cite{Bri1,Bri2} use a related
methodology to construct a scalar diffusion whose
one-dimensional marginal distributions are given
as a mixture of known densities.
Bentata and Cont~\cite{Ben}
recently announced an extension of
Gy\"ongy's result to jump diffusions under a
continuity assumption on the coefficients in the mimicking process and
a nondegeneracy assumption on the covariance or the jump measure of
the mimicking process.

Here we extend Gy\"ongy~\cite{G} in two ways.
First, we remove
the conditions of nondegeneracy and boundedness
on the covariance
of the It\^o process to be mimicked,
requiring only integrability
of this process and thereby extending the result to cover
popular stochastic volatility models such as the one due to
Heston~\cite{H}. Second,
we show that the mimicking process can preserve
the joint distribution of certain functionals of the It\^o
process (e.g., running maximum and running average)
at each fixed time. Our mimicking
process is a weak solution to an SDE, and in the case
of preservation of the joint distribution of functionals
of the It\^o process, the coefficients in this SDE may
depend on the values of
these functionals as well as the current value
of the underlying It\^o process.

The conditions that permit our construction are so weak that the
solution to the SDE we derive is not necessarily unique. Uniqueness
results, such as those found in Stroock and Varadhan~\cite{SV,SVV2},
require the conditional expectations determined by the Gy\"ongy--Dupire
formulas [see (\ref{GD}) in this paper] to be sufficiently regular
functions of the conditioning variables. It is difficult to see what
conditions one should impose on the data of our model (the processes
$b$ and $\sigma$ and the updating function $\Phi$ of Theorem
\ref{T2.4}) to ensure such regularity. Of course, if one
is willing to assume that the coefficients in the mimicking equation
are sufficiently well-behaved, then it is often possible to conclude
that the solution to the mimicking equation is unique.

Finally, we mention an independent body of work
devoted to a problem similar to the one considered here.
If an It\^o process is a submartingale,
Kellerer~\cite{Ke} has shown that it can be mimicked by
a Markov process. More generally,~\cite{Ke}
shows that given any set of marginal densities
$p(t,\cdot)$, $t\geq0$, that have finite first moments
and satisfy
$\int\varphi(y)p(s,y) \,dy\leq\int\varphi(y)p(t,y)$
for every $t\geq s\geq0$ and every nondecreasing
convex function $\varphi$, there is a Markov submartingale
whose density at each time $t$ is $p(t,\cdot)$.
Madan and Yor~\cite{MY} provide constructions of
such Markov processes in three specific cases in which
the first moments of $p(t,\cdot)$ are independent of $t$.
Cox, Hobson and Obl\'oj~\cite{Cox} and Ekstr{\"o}m et al.~\cite{Eks}
provide related constructions.
Forde~\cite{For} studies the problem of matching the joint law
of a process and its running maximum at an independent exponential time.
Our results address
the specific case in which
the densities $p(t,\cdot)$ are the marginals
of an It\^o process. Our mimicking
process satisfies an SDE,
but because the solution to this equation might not be unique,
we are not able to establish the Markov property
in all cases. On the other hand, we have the Gy\"ongy--Dupire
formulas for the drift and diffusion coefficients
of our mimicking process.

This paper is based on the first author's Ph.D. dissertation
\cite{Bru}. It is organized as follows. Section~\ref{DiscreteEx}
presents an intuitive discrete-time example that illustrates the main
ideas of our construction. In Section~\ref{SecMainResult} we state our
main result, Theorem~\ref{T2.4}, and provide some useful corollaries.
To prove Theorem~\ref{T2.4}, we construct a weakly relatively compact
sequence of processes that mimic some initial target process. We then
extract a limit from this sequence, check that the mimicking property
is preserved under weak convergence, and compute the semimartingale
characteristics of the limiting process. The tools to implement this
strategy are developed in Sections~\ref{SConMea}--\ref{S5}, and the
proof of Theorem~\ref{T2.4} is given in Section~\ref{S6}.

More specifically, in Section~\ref{SConMea} we begin with a probability
measure on path space and construct a ``concatenated'' measure which
assigns the same unconditional distribution as the original measure to
fragments of paths between concatenation time points but changes the
dependency structure across these time points. The new dependency
structure corresponds to ``partially forgetting'' the past at each
concatenation time point, and the resulting process possesses a limited
Markov-like property. The existence and uniqueness of the concatenated
measure are provided by Theorem~\ref{T3.3}, and Section~\ref{SExisUniq}
is devoted to the statement and proof of that theorem. Although the
concatenated measure may not be equivalent to the original measure,
certain properties of the process, such as finite variation and
absolute continuity, are preserved by the construction. The properties
we need are set out in Section~\ref{SPropPres}. The most important
result of this subsection is Proposition~\ref{P3.12}, which provides
conditions that are sufficient to ensure that the semimartingale
characteristics of the initial process are not disturbed by the
concatenation procedure.

Section~\ref{SecCondExp} sets out conditions under which the
conditional expectation of one process conditioned on a second process
can be written as a function of time and the second process. This
result, Proposition~\ref{P4.1}, is extended to include conditioning on
a random time as well in Proposition~\ref{P4.5}.

Finally, in Section~\ref{S5}, we set up the machinery for taking the
limit of a sequence of concatenated measures. Proposition~\ref{P5.2}
provides conditions on a sequence of weakly converging processes that
guarantee joint convergence of the processes and the integral of a
function of the processes. Proposition~\ref{P5.5} shows how to
approximate a process in $L_1$ by a piecewise constant process
constructed from the original process by sampling at random times.
Proposition~\ref{P5.7} shows that if a sequence of discrete-time
martingales is constructed by integrating with respect to time a
sequence of uniformly integrable processes and sampling these integrals
at stopping times, and if the maximum time between successive stopping
times approaches zero, then the integrand processes must also approach
zero.

\section{Guiding example} \label{DiscreteEx}

To motivate
the results that follow, we first
sketch a mimicking result for
discrete-time processes. This setting
illustrates the main ideas of our
proof methodology without the
technical complications of
continuous time.

Let $\N_0$ denote the set of
nonnegative integers,
let
$\sB(\R)$ denote
the Borel $\sigma$-field on $\R$
and let $(X_n)_{n\in\N_0}$ denote a
(not necessarily Markov)
stochastic process in discrete time that
takes values in $\mathbb R$. For each $n \in
\N_0$, we may
construct a measurable transition kernel $p_n\dvtx\mathbb
R \times\sB(\R) \rightarrow[0,1]$ with the property that $p_n(X_n;
A)$ is a version of $\P[X_{n+1} - X_n \in A | \sigma(X_n)]$ for
each $A
\in\sB(\R)$ and $A \mapsto p_n(x; A)$ is a probability measure for
each $x \in\mathbb R$.

After moving to a suitable extension of our
probability space if
necessary, we may construct a process
$Y$ such that $Y_0 = X_0$;
$Y_{n+1} - Y_n$ is conditionally
independent of $\sF_n$ given $Y_n$;
and $p_n(Y_n; A)$ is a version of
$\mathbb P(Y_{n+1} - Y_n \in A
| \sF_n)$ for each $n \in\mathbb N_0$
and $A \in\sB(\R)$. It
follows from these properties that
\[
\mathbb E\bigl[f(Y_{n+1}) | \sF_n\bigr] = \int
f(Y_n + x) p_n(Y_n; dx),
\]
so $Y$ is a Markov process. We also have
\[
\mathbb E\bigl[f(X_{n+1})\bigr] = \mathbb E \bigl[ \mathbb E\bigl[
f(X_{n+1}) | \sigma(X_n)\bigr] \bigr] = \mathbb E \biggl[
\int f(X_n + x) p_n(X_n; dx) \biggr]
\]
and $Y_0 = X_0$, so an inductive argument shows that $Y_n$ has the
same law as $X_n$ for each $n$. This is essentially the construction
given by Derman and Kani
\cite{DK2}.

Given a discrete-time process $X$, we
now let $\Xbar_n =
\max_{0\leq i\leq n} X_i$
denote the running maximum of the process $X$.
Although the law of the random variable $Y_n$ constructed above agrees
with the law of $X_n$ for each fixed $n$, the law of the process $Y$
may certainly differ from the law of the process~$X$. In particular,
the law of the pair $(X_n, \Xbar_n)$ may not agree with the law of the
pair $(Y_n, \Ybar_n)$ when $n \geq1$. Nevertheless, one can
construct a second process $Z$
such that the two-dimensional
process $(Z, \Zbar)$ is Markov and the joint law of the pair $(Z_n,
\Zbar_n)$ agrees with the joint law of the pair $(X_n, \Xbar_n)$ for
each~$n$, as we now show.

We
let $q_n\dvtx\R^2 \times\sB(\R) \rightarrow[0,1]$ denote a
transition kernel with the property that $p_n(X_n, \Xbar_n; A)$ is a
version of $\P[X_{n+1} - X_n \in A | \sigma(X_n, \Xbar_n)]$ for
each $A
\in\sB(\R)$. Moving to another extension of our probability space, we
may construct a process $Z$ such that $Z_0 = X_0$; $Z_{n+1} - Z_n$ is
conditionally independent of $\sF_n$ given $(Z_n, \Zbar_n)$; and
$p_n(Z_n, \Zbar_n; A)$ is a version of $\mathbb P(Z_{n+1} - Z_n \in A
| \sF_n)$
for each $n \in\mathbb N_0$.

We define
$\Phi\dvtx\mathbb{R}^3\rightarrow\mathbb{R}^2$
by $\Phi(e_1, e_2; x)
= (e_1+x, e_2 \vee(e_1 + x)) $, so that
$(Z_{n+1}, \Zbar_{n+1})
= \Phi(Z_n, \Zbar_n; Z_{n+1}-Z_n)$.
We may use the function $\Phi$
and the increments of the
process $Z$ to update the state of the process\vadjust{\goodbreak} $(Z, \Zbar)$. One
immediate consequence of this structure is that
\[
\mathbb E\bigl[f(Z_{n+1}, \Zbar_{n+1}) | \sF_n
\bigr] = \int f \circ\Phi(Z_{n}, \Zbar_{n}; y)
q_n(Z_{n}, \Zbar_{n}; dy),
\]
so $(Z, \Zbar)$ is a Markov process. We also have
\[
\mathbb E\bigl[f(X_{n+1}, \Xbar_{n+1})\bigr] = \mathbb E
\biggl[\int f \circ\Phi(X_{n}, \Xbar_{n}; y)
q_n(X_{n}, \Xbar_{n}; dy) \biggr],
\]
so another inductive argument shows that the law of the pair $(Z_n,
\Zbar_n)$ agrees with the law of the pair $(X_n, \Xbar_n)$ for each
$n$.
This paper
extends this
construction to continuous time.

\section{Main result}\label{SecMainResult}

In order to precisely state our main result,
we need some notation.
The symbol $\calE$ will always denote a
closed subset of a complete separable metric space,
that is, a Polish space.
Let $C^\calE$ be the space of
continuous functions from $[0,\infty)$
to $\calE$, endowed with the topology of
uniform convergence on compact
subsets of $[0,\infty)$.
We define the \textit{shift operator}
$\Theta\dvtx C^\calE\times\R\rightarrow
C^\calE$ by
\[
\Theta(x,t)\define x \bigl((t+\cdot)^+ \bigr),
\]
the \textit{stopping operator}
$\nabla\dvtx C^\calE\times[0,\infty)
\rightarrow C^\calE$ by
\[
\nabla(x,t)\define x(\cdot\wedge t)
\]
and, if $\sE$ is a vector space,
the \textit{difference operator}
$\Delta\dvtx C^\calE\times[0,\infty)
\rightarrow C^\calE$ by
\[
\Delta(x,t)\define x(t+\cdot)-x(t).
\]
In contrast to usual practice, here the shift operator can shift paths
to the right because $t$ can be negative, and in this case, the shifted
path takes the value $x(0)$ on $[0,-t]$. The difference operator
actually maps into $C_0^\calE$, the space of continuous functions from
$[0,\infty)$ to $\calE$ with initial condition zero. If $\calE=\R^d$
for some integer~$d$, we write $C^d$ and $C^d_0$ rather than $C^{\R^d}$
and $C^{\R^d}_0$.

Fix a Polish space $\calE$, fix a positive integer $d$ and define
$\Omega^{\calE,d}\define\calE\times C_0^d$. We endow $\Omega^{\calE,d}$
with the product topology. We denote a generic element of
$\Omega^{\calE,d}$ by $\omega=(e,x)$ and define the random variable
$E(e,x)=e$ and the $\R^d$-valued process $X(e,x)=x$. For a random time
$T$, we use the notation $X^T$ to denote the process $X$ stopped at
$T$, that is,
%
%
\begin{equation}
\label{2.4} X^T_t(\omega)=X_{t\wedge T(\omega)}(\omega) =
\nabla_t \bigl(X(\omega),T(\omega) \bigr),\qquad t\geq0.
\end{equation}

%
\begin{definition}\label{D2.1}
We say that
$\Phi\dvtx\Omega^{\calE,d}\rightarrow C^\calE$
is an \textit{updating
function} provided
%
%
\begin{eqnarray}
\label{2.5}
\Phi_0(e,x)&=&e,\qquad e\in\calE,
\\
\label{2.6}
\Phi^t(e,x)&=&\Phi^t \bigl(e,\nabla(x,t) \bigr),\qquad t\geq0,
e\in\calE, x\in C_0^d,
\\
\label{2.7}
\Theta\bigl(\Phi(e,x),t \bigr) &=&\Phi\bigl(\Phi_t(e,x),
\Delta(x,t) \bigr),\qquad t\geq0, e\in\calE, x\in C_0^d.
\end{eqnarray}
\end{definition}

In other words, $\Phi$ takes an initial
condition in $\calE$ [see (\ref{2.5})]
and a path in $C_0^d$ and
generates a path in $C^\calE$.
Property (\ref{2.6}) says that the path
$\Phi(e,x)$ stopped
at $t$ depends only on the initial condition $e$
and the path of $x$ stopped at $t$. This is
a nonanticipative
property. Property (\ref{2.7}) is a type
of Markov property,
but on a path-by-path basis without the
presence of a probability
measure. It implies that the path of
$\Phi(e,x)$ from time $t$
onward depends only on the value of the
path at time $t$
and the increments of $x$ from time $t$
onward. Using
the characterization of the Markov property
as independence
of the future and past given the present,
it is easily verified that
if $\xi$ is a continuous $\R^d$-valued Markov
process, and if
for each $t$ the value of $\xi_t$
can be deduced from the value
of $\Phi_t(\xi_0,\xi-\xi_0)$, then $\Phi(\xi_0,\xi-\xi_0)$ is
also Markov.\looseness=1

%
%
\begin{example}[(Process itself)]\label{Ex2.3}
A trivial case of an updating
function is obtained if we let
$\calE=\R^d$, $\Omega^{\calE,d}=\calE\times C_0^d$ and
$\Phi(e,x)=e+x$ for $e\in\R^d$ and $x\in C_0^d$.
If $\xi$ is a continuous
$\R^d$-valued Markov process and we
represent $\xi$ as
$(\xi_0,\xi-\xi_0)\in\calE\times C_0^d$,
then $\Phi_t(\xi_0,\xi-\xi_0)=\xi_t$
and $\Phi(\xi_0,\xi-\xi_0)=\xi$ is Markov.
\end{example}

%
%
\begin{example}[(Integral-to-date)]\label{IntegralEx}
Let $\calE=\R^2$ and $\Omega^{\calE,1}=\calE\times C_0^1$. We interpret
a point $(e_1,e_2; x)\in\Omega^{\sE,1}$ as a path $e_1+x$ with initial
condition $e_1+x(0)=e_1$ and the initial value of a running integral
given by $e_2$. It is then easy to check that
\[
\Phi_t(e_1,e_2; x) = \biggl(e_1+x(t),
e_2 + \int_0^t
\bigl(e_1+x(s) \bigr) \,ds \biggr)
\]
is an updating function.
\end{example}

%
%
\begin{example}[(Maximum-to-date)]\label{Ex2.2}
Let $\calE=\{(e_1,e_2)\in\R^2\dvtx e_1\leq e_2\}$ and
$\Omega^{\calE,1}=\calE\times C_0^1$. We regard the generic element
$(e_1,e_2; x)\in\Omega^{\sE,1}$ as a path $e_1+x$ with initial
condition $e_1+x(0)=e_1$ and the time-zero maximum-to-date $e_2$. Given
such a triple, the value of the path at a later time $t$ and the
maximum-to-date at that time $t$ are $e_1+x(t)$ and $e_2\vee\max_{0\leq
s\leq t} (e_1+x(s) )$, respectively. We thus define
\[
\Phi_t(e_1,e_2; x) = \Bigl(e_1+x(t),
e_2\vee\max_{0\leq s\leq t} \bigl(e_1+x(s) \bigr)
\Bigr).
\]
It is straightforward to verify that $\Phi$ is an updating function. If
$\xi$ is a continuous real-valued Markov process, then
$\Phi_t(\xi_0,M_0; \xi-\xi_0) =(\xi_t, M_0\vee\max_{0\leq s\leq
t}\xi_s)$ is also Markov, where $M_0$ is any random variable satisfying
$M_0\geq\xi_0$ almost surely.
\end{example}

As a final extremal example,
we give an updating function that
records the entire history of the path.
%
%
\begin{example}[(Path-to-date)]\label{Ex2.4}
Define $\calE= \{(s,x)\in[0,\infty)\times C^d; x$
is constant on $[s,\infty) \}$, define $\Omega^{\calE,d}=\calE\times C_0^d$, and set
\begin{eqnarray}
\Phi_t(s,x; y)= \bigl(s+t, \nabla\bigl(\nabla(x,s) +
\Theta(y,-s),s+t \bigr) \bigr),\nonumber\\
&&\eqntext{x\in C^d, s\in[0,\infty), y\in
C_0^d.}
\end{eqnarray}
Given paths $x\in C^d$ and $y\in C_0^d$
and a time $s\geq0$,
$\nabla(x,s)+\Theta(y,-s)$
is the path that follows $x$ on
$[0,s]$ with
$y$ appended after time $s$. The second component
of $\Phi_t$
is this path stopped at time $s+t$.
The first component
of $\Phi_t$ is the time $s+t$ at which
this path is stopped.
As $t$ marches forward, the second component
of the
operator $\Phi$ applied to $(s,x; y)$
appends more and more of the path $y$ to
the path $x$,
always appending at time $s$. It is tedious
but straightforward to check that $\Phi$
is an updating function. For any continuous
$\R^d$-valued process $\xi$, we have
%
%
\begin{equation}
\label{2.12a} \Phi_t\bigl(0,\xi^0; \xi-\xi_0
\bigr)=\bigl(t,\xi^t\bigr),\qquad t\geq0,
\end{equation}
where we recall from (\ref{2.4}) that
$\xi^t$ is the process
$\xi$ stopped at $t$.
\end{example}

%
%
\begin{theorem}[(Main result)]\label{T2.4} Suppose an $\R^d$-valued
process $Y$ is~given~by
%
%
\begin{equation}
\label{2.10} Y_t=\int_0^tb_s
\,ds + \int_0^t\sigma_s
\,dW_s,\qquad t\geq0,
\end{equation}
where $W$ is an $\R^r$-valued Brownian motion under some probability
measure $\P$, $b$ is an $\R^d$-valued process adapted to a filtration
under which $W$ is a Brownian motion and $\sigma$ is a $d\times r$
matrix-valued process adapted to the same filtration as $b$. Let
$\calE$ be a Polish space, define $\Omega^{\calE,d}\define\calE\times
C_0^d$, let $\Phi\dvtx\Omega^{\calE,d}\rightarrow C^{\calE}$ be a
continuous updating function, let $Z_0$ be an $\calE$-valued random
variable and set $Z=\Phi(Z_0,Y)$, which is a continuous $\calE$-valued
process. Finally,\vspace*{1pt} assume that $\E\int_0^t (\|b_s\|
+\|\sigma_s\sigma_s^{tr}\| ) \,ds<\infty$ for $t\geq0$. Then there
exists an $\R^d$-valued measurable function $\widehat{b}$ and a
$d\times d$ matrix-valued measurable function $\sigmahat$, both defined
on $[0,\infty)\times\calE$, and there exists a Lebesgue-null
set $N\subset[0,\infty)$, so that\setcounter{footnote}{2}\footnote{We interpret (\ref{GD}) and subsequent similar
equations to mean
that for each fixed $t\in N^c$, the left-hand
side of each equation is a version
of the conditional expectation appearing on the right-hand
side.}%
%
%
\begin{eqnarray}
\label{GD} \widehat{b}(t,Z_t)&=& \E[b_t | Z_t],\nonumber\\[-8pt]\\[-8pt]
\sigmahat(t,Z_t)\sigmahat^{tr}(t,Z_t)&=& \E\bigl[
\sigma_t\sigma_t^{tr} | Z_t\bigr],\qquad
\P\mbox{-a.s.}, t\in N^c.\nonumber
\end{eqnarray}
Furthermore,\vspace*{2pt} there exists a filtered probability space $(\Omegahat,
\sFhat,\{\sFhat_t\}_{t\geq0},\Phat)$ that supports a continuous
$\R^d$-valued adapted process $\Yhat$, a continuous $\calE$-valued
adapted process $\Zhat$ and a $d$-dimensional Brownian motion $\What$
satisfying
%
%
\begin{equation}\label{2.14}
\Yhat_t = \int_0^t\widehat{b}(s,
\Zhat_s) \,ds +\int_0^t\sigmahat(s,
\Zhat_s) \,d\What_s,\qquad \Zhat= \Phi(\Zhat_0,
\Yhat),\qquad t\geq0,\hspace*{-28pt}
\end{equation}
and such that for each $t\geq0$, the distribution of $\Zhat_t$ under
$\Phat$ agrees with the distribution of $Z_t$ under $\P$.
\end{theorem}

Although both $Y$ in (\ref{2.10})
and $\Yhat$ in (\ref{2.14}) are $d$-dimensional
processes, the ``state'' $\Zhat$ of the system
in (\ref{2.14}) can be of a much lower
dimension than the state process
needed to describe (\ref{2.10}).
In (\ref{2.10})
the processes $b$ and $\sigma$ are typically
given by stochastic differential equations
driven by additional Brownian motions not
mentioned in the statement of the theorem.
The process $\Zhat$ is typically the
process $\Yhat$ itself augmented by
some functional of the path of $\Yhat$.
We give examples below. Indeed, the
remainder of this section
illustrates the applications of Theorem
\ref{T2.4}.
In this section we also show
by example that (\ref{2.14}) can have multiple
solutions
and discuss
conditions that guarantee uniqueness.
The subsequent sections are
devoted to the proof of
Theorem~\ref{T2.4}.

Note\vspace*{2pt} that $Y$ in Theorem~\ref{T2.4} is a martingale
if and only if $b_s$ is zero for Lebesgue almost every
$s$ almost surely. In this case,\vspace*{1pt} $\widehat{b}$
is also zero, and $\widehat{Y}$ is a local martingale.
But since $\widehat{Z}_s$ has the same distribution
as $Z_s$ for each $s$, the integrabilty condition
assumed on $\sigma\sigma^{tr}$ implies the same condition
on $\widehat{\sigma}\widehat{\sigma}^{tr}$
and $\widehat{Y}$ is in fact a martingale.

As a first application, we take $Y = X - X_0$ and $Z=X$ in Theorem
\ref{T2.4} and use the updating function of Example~\ref{Ex2.3}. We
then have the following corollary, which
is the result obtained by
Gy\"{o}ngy~\cite{G}, but here without the
boundedness and
nondegeneracy
assumptions of~\cite{G}.

%
%
\begin{corollary}[(Process itself)]\label{C2.6}
Suppose an $\R^d$-valued process $X$ is given by
%
%
\begin{equation}
\label{2.11gb} X_t = X_0 + \int_0^tb_s
\,ds +\int_0^t\sigma_s
\,dW_s,\qquad t\geq0,
\end{equation}
where $W$, $\P$, $b$ and $\sigma$ are as in Theorem~\ref{T2.4}. Then
there exists an $\R^d$-valued measurable function $\widehat{b}$ and
a $d\times d$ matrix-valued measurable function $\sigmahat$,
both defined on $[0,\infty)\times\R^d$, and there exists a
Lebesgue-null set $N$, so that
%
%
\begin{eqnarray}
\label{2.18} \widehat{b}(t,X_t)&=& \E[b_t |
X_t],\nonumber\\[-8pt]\\[-8pt] \sigmahat(t,X_t)\sigmahat^{tr}(t,X_t)&=&
\E\bigl[\sigma_t\sigma_t^{tr} |
X_t \bigr],\qquad \P\mbox{-a.s.}, t\in N^c.\nonumber
\end{eqnarray}
Furthermore,
there exists a filtered probability space
$(\Omegahat,\sFhat,\{\sFhat_t\}_{t\geq0},\Phat)$
that supports a continuous
$\R^d$-valued adapted process $\Xhat$ and
a $d$-dimensional Brownian motion
$\What$ satisfying
%
%
\begin{equation}
\label{2.19} \Xhat_t=\Xhat_0+\int_0^t
\widehat{b}(s,\Xhat_s) \,ds +\int_0^t
\sigmahat(s,\Xhat_s) \,d\What_s,\qquad t\geq0,
\end{equation}
and such that for each $t\geq0$,
the distribution of $\Xhat_t$
under $\Phat$ agrees with the distribution
of $X_t$ under $\P$.
\end{corollary}

%
%
\begin{example}[(Fake Brownian motion)]
Let $G_1$ and $G_2$ be standard normal random variables, let
$(B_t)_{t \geq0}$ be Brownian motion and assume that
$G_1$, $G_2$, and $W$ are independent. Define the process
\[
X_t = \sqrt{t} (G_1 \cos B_{\ln t} +
G_2 \sin B_{\ln t}),\qquad t\geq1,
\]
and set $\sF_t = \sigma(G_1, G_2, B_{s}, 0 \leq s \leq\ln t)$ for
$t \geq1$. Then $X$ is a continuous martingale with respect to
$\{\sF_t\}_{t \geq1}$ and $\langle X \rangle_t = \int_1^t
\sigma^2_s \,ds$ for $t \geq1$, where $\sigma_t = -G_1 \sin
B_{\ln t} + G_2 \cos B_{\ln t}$. In particular, we may write $X$
in the form
\[
X_t = X_1 + \int_1^t
\sigma_s \,dW_s,\qquad t \geq1,
\]
for some Brownian motion $(W_t)_{t\geq1}$.

Conditioned on the value of $B_{\ln t}$, the
random variables $X_t /\sqrt{t}$ and $\sigma_t$ are
independent and
standard normal, so they are
unconditionally independent and standard normal.
Consequently,
\[
\E\bigl[\sigma^2_t | X_t\bigr] = \E\bigl[
\sigma^2_t\bigr] = 1,\qquad \P\mbox{-a.s.}, t \geq1,
\]
and we may take $\widehat{b}= 0$ and $\sigmahat= 1$ in the previous
corollary. As $X_1$ is standard normal, the previous corollary,
adapted to the time interval $[1,\infty)$, asserts that the process
$X$ has the same one-dimensional marginal distributions as a
Brownian motion on $[1,\infty)$. This is
not hard to check directly in this example.

This construction is due to Oleszkiewicz~\cite{O} who was interested in
producing a \textit{fake Brownian motion} (see also~\cite{Al,Ha}). A
fake Brownian motion is a continuous martingale that has the same
one-dimensional marginal distributions as a Brownian motion but is not
itself a Brownian motion. Oleszkiewicz shows that the process $X$
constructed above can be extended to produce a fake Brownian motion on
the time interval $[0,\infty)$. The argument given in this example
can be extended to show that the process which mimics Oleszkiewicz's
fake Brownian motion in the sense of Corollary~\ref{C2.6} is simply
Brownian motion.
\end{example}

Taking $Y = X - X_0$ and $Z_t = (X_t, A_t)$,
and using the
updating function in Example~\ref{IntegralEx}, we
obtain the following corollary about the
distribution of a process and its running integral.

%
%
\begin{corollary}[(Integral-to-date)]\label{IntegralCor}
Suppose a real-valued process $X$ is given by
(\ref{2.11gb})
where $W$, $\P$, $b$ and $\sigma$
are as in Theorem~\ref{T2.4}\vadjust{\goodbreak} with $d=r=1$.
Let $A$ be a continuous process such that
\[
A_t=A_0 + \int_0^t
X_s \,ds,\qquad t\geq0.
\]
Then there exists a real-valued measurable
function $\widehat{b}$ and a
$[0,\infty)$-valued
measurable function $\sigmahat$,
both defined on
$[0,\infty)\times\mathbb R^2$,
and there exists a Lebesgue-null
set $N$, such that
\begin{eqnarray*}
\widehat{b}(t,X_t,A_t)&=& \E[b_t |
X_t,A_t],\\ \sigmahat^2(t,X_t,A_t)&=&
\E\bigl[\sigma_t^2 | X_t,A_t
\bigr],\qquad \P\mbox{-a.s.}, t\in N^c.
\end{eqnarray*}
Furthermore,
there exists a filtered probability space
$(\Omegahat,\sFhat,\{\sFhat_t\}_{t\geq0},
\Phat)$ that supports continuous
real-valued adapted processes $\Xhat$ and $\Ahat$
and a real-valued Brownian motion
$\What$ satisfying
%
%
\begin{eqnarray}
\label{eq:integraltodate}\qquad
\Xhat_t &=&
\Xhat_0+\int_0^t\widehat{b}(s,
\Xhat_s, \Ahat_s) \,ds +\int_0^t
\sigmahat(s,\Xhat_s,\Ahat_s) \,d\What_s,\qquad t
\geq0,
\nonumber\\[-8pt]\\[-8pt]
\Ahat_t &=& \Ahat_0 + \int_0^t
\Xhat_s \,ds,\qquad t\geq0, \nonumber
\end{eqnarray}
and such that for each $t\geq0$,
the distribution of the pair $(\Xhat_t, \Ahat_t)$
under $\Phat$ agrees with the distribution
of the pair $(X_t, A_t)$ under $\P$.
\end{corollary}

Taking $Y = X - X_0$ and $Z_t = (X_t, M_t)$,
and using the
updating function in Example~\ref{Ex2.2}, we
obtain the following corollary about the
distribution of a process and its running maximum.

%
%
\begin{corollary}[(Maximum-to-date)]\label{C2.7}
Suppose a real-valued process $X$ is given by
(\ref{2.11gb})
where $W$, $\P$, $b$ and $\sigma$
are as in Theorem~\ref{T2.4} with $d=r=1$.
Let $M_0$ be a random variable
satisfying $M_0\geq X_0$ almost surely
and define
\[
M_t=M_0\vee\max_{0\leq s\leq t}X_s,\qquad t\geq0.
\]
Then there exists a
real-valued measurable
function $\widehat{b}$ and a
$[0,\infty)$-valued
measurable function $\sigmahat$,
both defined on
$[0,\infty)\times\mathbb R^2$,
and there exists a Lebesgue-null
set $N$, such that
\begin{eqnarray*}
\widehat{b}(t,X_t,M_t)&=& \E[b_t |
X_t,M_t],\\
\sigmahat^2(t,X_t,M_t)&=&
\E\bigl[\sigma_t^2 | X_t,M_t
\bigr],\qquad \P\mbox{-a.s.}, t\in N^c.
\end{eqnarray*}
Furthermore,
there exists a filtered probability space
$(\Omegahat,\sFhat,\{\sFhat_t\}_{t\geq0},\Phat)$
that supports continuous
real-valued adapted processes $\Xhat$ and $\Mhat$
and a real-valued Brownian motion
$\What$ satisfying
%
%
\begin{eqnarray}
\label{eq:runningmax} \qquad
\Xhat_t &=&
\Xhat_0+\int_0^t\widehat{b}(s,
\Xhat_s, \Mhat_s) \,ds +\int_0^t
\sigmahat(s,\Xhat_s,\Mhat_s) \,d\What_s,\qquad t
\geq0,
\nonumber\\[-8pt]\\[-8pt]
\Mhat_t &=& \Mhat_0\vee\max_{0\leq s\leq t}
\Xhat_s,\qquad t\geq0,\nonumber
\end{eqnarray}
and such that for each $t\geq0$,
the distribution of the pair $(\Xhat_t, \Mhat_t)$
under $\Phat$ agrees with the distribution
of the pair $(X_t, M_t)$ under $\P$.
\end{corollary}

Taking $Y = X - X_0$ and $Z_t = (t,X^t)$,
and using the
updating function in Example~\ref{Ex2.4}, we
obtain the following corollary, which states that
every It\^o process with integrable drift and
covariance is a weak solution to an
SDE with path-dependent coefficients.

%
%
\begin{corollary}[(Path-to-date)]
Suppose a real-valued process $X$ is given by
(\ref{2.11gb})
where $W$, $\P$, $b$ and $\sigma$
are as in Theorem~\ref{T2.4}.
Then there exist
path-dependent functionals $\widehat{b}$ and
$\sigmahat$, both defined on $[0,\infty) \times
C^d$, with $\widehat{b}$ taking values in $\R^d$ and
$\sigmahat$ taking values in the space of $d\times
d$ matrices and a Lebesgue-null set $N$ such that
\begin{eqnarray*}
\widehat{b}\bigl(t,X^t\bigr)&=& \E\bigl[b_t |
X^t\bigr],\\
\sigmahat\bigl(t,X^t\bigr)
\sigmahat^{tr}\bigl(t,X^t\bigr)&=& \E\bigl[
\sigma_t\sigma_t^{tr} | X^t\bigr],\qquad
\P\mbox{-a.s.}, t\in N^c.
\end{eqnarray*}
Furthermore,
there exists a filtered probability space
$(\Omegahat,\sFhat,\{\sFhat_t\}_{t\geq0},\Phat)$
that supports a continuous
$\R^d$-valued adapted process $\Xhat$ and
a $d$-dimensional Brownian motion
$\What$ satisfying
\[
\Xhat_t=\Xhat_0+\int_0^t
\widehat{b}\bigl(s,\Xhat^s\bigr) \,ds +\int_0^t
\sigmahat\bigl(s,\Xhat^s\bigr) \,d\What_s,\qquad t\geq0,
\]
and such that $\Xhat$ has the same distribution
under $\Phat$ as $X$ has under $\P$.
\end{corollary}

We close this section with a brief discussion of the nonuniqueness
that can arise in equation (\ref{2.14}) of Theorem~\ref{T2.4} and its
relationship to the strong Markov property. We first provide a simple
example within the context of Corollary~\ref{C2.6}, where $X$ and
$Z$ are the same process.

%
%
\begin{example}[(Nonuniqueness)]\label{Ex2.8}
Let $d=1$ and $b=0$ in Corollary~\ref{C2.6} and let
$X_t=\int_0^t\sigma_s \,dW_s$,
where
$\sigma_s
=\ind_{(1,\infty)}(s)\ind_{\{W_1>0\}}$.
Then
$X_t=\ind_{(1,\infty)}(t)\*
\ind_{\{W_1>0\}}(W_t-W_1)$.
From (\ref{2.18}) we see that
$\sigmahat(t,y)=0$
for $0\leq t\leq1$, and for $t>1$,
\[
\sigmahat^2(t,y)=\E\bigl[\sigma_t^2 |
X_t=y\bigr] =\cases{ %
1,&\quad if $y\neq0$,
\cr
0,&\quad if $y=0$.}
\]
Both $\Xhat^1_t\equiv0$ and
$\Xhat^2_t=\ind_{(1,\infty)}(t)(W_t-W_1)$
are solutions of (\ref{2.19}).
The weak solution $\Xhat$ that has
the same one-dimensional distributions
as $X$ is obtained by an initial randomization
that is independent of $W$
and determines whether $\Xhat$
agrees with $\Xhat^1$ or\vadjust{\goodbreak}
$\Xhat^2$, each of these
events having probability~$\frac12$.
This process is Markov, but not strong
Markov, as can be seen by considering the stopping
time that is the first time after time $2$
that zero is reached.
\end{example}

The previous example shows that the mimicking
process may not be strong
Markov. Nevertheless, if we are willing to
impose further conditions
on the coefficients $\widehat{b}$ and $\sigmahat$
appearing in
Theorem~\ref{T2.4}, then we can often
conclude that
the solution to (\ref{2.19}) is unique in law and
strong Markov.
In particular,
if we assume that
$\widehat{b}$ appearing
in Corollary~\ref{C2.6} is
bounded and measurable and
that
$\sigmahat\sigmahat^{tr}$ is bounded,
strictly positive-definite and continuous, then
the results of Stroock and Varadhan~\cite{SV,SVV2}
ensure that the
mimicking process satisfying (\ref{2.19})
in Corollary~\ref{C2.6} is
unique in law and
strong Markov with respect
to its natural filtration. We state this
observation as a corollary.
%
%
\begin{corollary}
Let $X$ denote an $\R^d$-valued process that satisfies equation
(\ref{2.11gb}), where $W$, $\P$, $b$ and $\sigma$ are as in Theorem
\ref{T2.4}, and suppose that there exists a locally bounded
measurable $\R^d$-valued function $\widehat{b}$ and a measurable
$d\times
d$ matrix-valued function $\sigmahat$ such that (\ref{2.18}) holds
and the function $\widehat a(t,x) = \sigmahat(t,x)\sigmahat^{tr}(t,x)$ is
continuous and strictly positive definite.
Then there exists
a weak solution to the SDE (\ref{2.19}) and all weak solutions
have the same law. Moreover, if
$\Xhat$ is a weak solution to (\ref{2.19}), then $\Xhat$ is strong
Markov with respect to the filtration $\sFhat_t =\sigma(\Xhat_s, 0
\leq s \leq t)$ and has the same one-dimensional marginal
distributions as the process $X$.
\end{corollary}

The conditions in this corollary can be
weakened. For example, more recent
results of Krylov~\cite{Kry4} imply
that the mimicking process
in Corollary~\ref{C2.6} is
unique in law and strong Markov when $\widehat{b}$ is
bounded and measurable and
$\sigmahat\sigmahat^{tr}$ is bounded,
locally uniformly positive-definite and
continuous in the sense of
vanishing mean oscillation. If we
restrict attention
to the one-dimensional case, then the
mimicking process
in Corollary~\ref{C2.6} is
unique in law and strong Markov when $\widehat{b}$ is
bounded and measurable, and $\sigmahat$ is bounded,
locally uniformly
positive and measurable
(Exercise 7.3.3 of~\cite{SV2}).

The two-dimensional
process $(\Xhat, \Ahat)$ in Corollary~\ref{IntegralCor}
is degenerate, so the results of Stroock
and Varadhan~\cite{SV,SVV2,SV2} do not apply.
However, Theorem~5.10 of
\cite{Bru2} asserts that
the solution to (\ref{eq:integraltodate})
is uniquely determined in law
when $\widehat{b}$ is
bounded and measurable and $\sigmahat$ is bounded,
strictly positive and
continuous. It then follows under these conditions
that the pair of mimicking processes
in Corollary~\ref{IntegralCor} possess the strong
Markov property.

We observe finally that the path
functional
$x \mapsto\max_{s \in[0,t]} x(s)$ is Lipschitz
continuous for each fixed $t\geq0$. This implies
that pathwise
uniqueness holds for the mimicking
equation (\ref{eq:runningmax})
in Corollary~\ref{C2.7} when $\widehat{b}$ and $\sigmahat$ are
bounded and locally Lipschitz continuous.\vadjust{\goodbreak}
As a result, it is easy to
check that the process $(\Xhat, \Mhat)$ in
Corollary~\ref{C2.7} is strong Markov
under these conditions.

To summarize, we cannot conclude
in general that the mimicking process
$\Zhat$
in Theorem~\ref{T2.4} is unique in law and
strong Markov.
In many cases
of interest, it is
possible to identify conditions
that may be imposed on the mimicking
equation to ensure that the solution
is unique and that the mimicking
process possesses the strong Markov property.
However, these
conditions vary from case to case,
and depend in an essential way on
the structure of the updating function.

\section{Concatenated measure}\label{SConMea}

In this section we begin with
a measure $\P$
and a partition $\Pi$ of
$[0,\infty)$ and construct a
\textit{concatenated measure}.
This is the continuous-time analogue
of the measure induced on path space
by the process $Y$ or the pair
$(Z, \overline{Z})$ in Section~\ref{DiscreteEx}.
We use the notation introduced
at the beginning of
Section~\ref{SecMainResult}.
On the space $\Omega^{\calE,d}=\calE\times C_0^d$,
we introduce
the $\sigma$-field
${\cal F}^{\sE,d}\define\frakE\otimes\sigma(X)$ and the
filtration
${\cal F}^{\sE,d}_t\define\frakE\otimes\sigma(X^t)$,
$t\geq0$,
where $\frakE$ is the Borel
$\sigma$-field in $\calE$.

%
%
\begin{definition}\label{D3.2}
Let $0=T_0\leq T_1\leq\cdots\leq T_n$
be a sequence of
finite (for \textit{every} $\omega$)
$\{{\cal F}^{\sE,d}_t\}_{t\geq0}$-stopping times and let
$\{\sG_i\}_{i=0}^n$ be a collection of
$\sigma$-fields satisfying
$\sG_i\subset{\cal F}^{\sE,d}_{T_i}$ for
$i=0,\ldots,n$. Set $T_{n+1}=\infty$,
set $\sH_0={\cal F}^{\sE,d}_0$ and define
$\sH_{i+1}\define
\sG_i\vee\sigma
(\Delta(X^{T_{i+1}},T_i) )$,
$i=0,1,\ldots,n$.
We say that
$\Pi\define(T_i,\sG_i)_{i=0}^n$
is an \textit{extended partition} provided:
\begin{longlist}[(a)]
\item[(a)]
$T_{i+1}-T_i\in\sG_i\vee
\sigma(\Delta(X,T_i) ),
i=0,1,\ldots,n-1$,
\item[(b)]
$\sG_i\subset\sH_i, i=0,1,\ldots,n$.
\end{longlist}
\end{definition}

%
%
\begin{remark}\label{R3.2}
Because $T_{i+1}-T_i$ is
${\cal F}^{\sE,d}_{T_{i+1}}$-measurable
and ${\cal F}^{\sE,d}_{T_{i+1}}=\frakE\otimes
\sigma(X^{T_{i+1}})$,
condition (a) in Definition~\ref{D2.1}
is equivalent to the
apparently stronger condition:
\begin{longlist}[(a$'$)]
\item[(a$'$)]
$T_{i+1}-T_i\in\sH_{i+1}, i=0,1,\ldots,n-1$.
\end{longlist}
Because $\sG_i\subset{\cal F}^{\sE,d}_{T_i}
\subset{\cal F}^{\sE,d}_{T_{i+1}}$
and $\sigma(\Delta(X^{T_{i+1}},T_i))
\subset{\cal F}^{\sE,d}_{T_{i+1}}$,
we have $\sH_{i+1}\subset{\cal F}^{\sE,d}_{T_{i+1}}$, or equivalently,
%
%
\begin{equation}
\label{3.1a} \sH_{i}\subset{\cal F}^{\sE,d}_{T_{i}},\qquad
i=0,1,\ldots,n,n+1.
\end{equation}
\end{remark}

An extended partition is a model
for observing and partially forgetting
information over time. Partial forgetting
occurs in Section~\ref{DiscreteEx} when
we condition on the value of a process
at time $n$ rather than on $\sF_n$.
With an extended partition, at time $T_i$
we retain the information in $\sG_i$ as
we move forward into the interval $[T_i,T_{i+1}]$,
but carry no other
information from ${\cal F}^{\sE,d}_{T_i}$ forward. We
then observe increments in
$X$ over the interval $[T_i,T_{i+1}]$,
so that the information we have
at time $T_{i+1}$ is $\sH_{i+1}$. This
information is sufficient to
tell us the length of time $T_{i+1}-T_i$
we conduct the observations.
We then remember
only the information in the
sub-$\sigma$-field $\sG_{i+1}$ of
$\sH_{i+1}$ as we go forward into the
interval $[T_{i+1},T_{i+2}]$.

\subsection{Existence and uniqueness
of concatenated measure}\label{SExisUniq}

%
%
\begin{theorem}[($n$-fold concatenation)] \label{T3.3}
Let $\P$ be a probability measure on
$(\Omega^{\calE,d},{\cal F}^{\sE,d})$,
and let $(T_i,\sG_i)_{i=0}^n$
be an extended partition.
Then there exists a unique measure
$\P^{\otimes\Pi}$ satisfying
%
%
\begin{eqnarray}
\label{3.2}
\P^{\otimes\Pi}[A]&=&\P[A],\qquad A\in\sH_i,\qquad i=0,1,\ldots,n+1,
\\
\label{3.3}
\P^{\otimes\Pi} \bigl[B | {\cal F}^{\sE,d}_{T_i} \bigr] &=&
\P[B | \sG_i ],\qquad B\in\sH_{i+1},\qquad i=0,1,\ldots,n.
\end{eqnarray}
We interpret (\ref{3.3}) to mean that every $\P$-version of $\P[B |
\sG_i]$ is a $\P^{\otimes\Pi}$-version of $\P^{\otimes\Pi}[B | {\cal
F}^{\sE,d}_{T_i}]$.
\end{theorem}

%
%
\begin{example}[(Simple concatenated
measure)]\label{SimpExPart1} Let $\sE=\{0\}$, so that $\Omega^{\sE,1}$
is isomorphic to $C_0^1$. Then $\sF^{\sE,1}_0$ is the trivial
$\sigma$-algebra $\{\varnothing,C_0^1\}$. We consider the extended
partition $\Pi=(T_i,\sG_i)_{i=0}^{1}$ with
$\sG_0=\sG_1=\{\varnothing,C_0^1\}$ and $T_0=0$, $T_1=1$ and, by
convention, $T_2=\infty$. Then $\sH_0=\sF^{\sE,1}_0$,
$\sH_1=\sF^{\sE,1}_1=\sigma(X(t), 0\leq t\leq1)$ and
$\sH_2=\sigma(X(t)-X(1), t\geq1)$. We define four elements of $C_0^1$
by $\omega^0(t)=0$, $\omega^1(t)=t$, $\omega^2(t)=t\wedge1$ and
$\omega^3(t)=(t-1)^+$ for $t\geq0$. Let $\delta^i$ be the probability
measure on $C_0^1$ assigning probability $1$ to $\omega^i$, and set
$\P=(\delta^0+\delta^1)/2$. The sets
\[
A_0= \bigl\{x\in C_0^1\dvtx x(t)=0\ \forall t
\in[0,1] \bigr\} \quad\mbox{and}\quad A_1= \bigl\{x\in C_0^1\dvtx x(t)=t
\ \forall t\in[0,1]\bigr\}
\]
are in $\sH_1=\sF^{\sE,1}_1$
and $\P(A_0)=\P(A_1)=\frac12$.
According to (\ref{3.2}), we must also have
$\P^{\otimes\Pi}(A_0)=\P^{\otimes\Pi}(A_1)=\frac12$.
The sets
\begin{eqnarray*}
B_0 &=& \bigl\{x\in C_0^1\dvtx x(t)-x(1)=0
\ \forall t\in[1,\infty)\bigr\},
\\
B_1 &=& \bigl\{x\in C_0^1\dvtx x(t)-x(1)=t-1
\ \forall t\in[1,\infty)\bigr\}
\end{eqnarray*}
are in $\sH_2$, and
(\ref{3.3}) implies that
\[
\P^{\otimes\Pi}\bigl[B_0 | \sF^{\sE,1}_1\bigr]
=\P[B_0 | \sG_1] =\P[B_0] =\tfrac12.
\]
Integrating this equation over $A_1$
with respect to $\P^{\otimes\Pi}$,
we see that
\[
\tfrac14=\P^{\otimes\Pi}(A_1\cap B_0)=\P^{\otimes\Pi}
\bigl(\omega^2\bigr).
\]
Considering all combinations of $A_j$
and $B_k$, we conclude that
$\P^{\otimes\Pi}(\omega^i)=\frac14$ for $i=0,1,2,3$,
that is, $\P^{\otimes\Pi}=(\delta_0+\delta_1+\delta_2+\delta_3)/4$.
\end{example}

The remainder of this subsection is devoted
to the proof of Theorem~\ref{T3.3}. Let $\frakC^d$ denote the Borel
$\sigma$-field in $C^d$ and let $\frakC_0^d$ denote the trace
$\sigma$-field in $C_0^d$. We first concatenate a deterministic
initial path and a probability measure at a deterministic time. Given
a fixed point $\omegabar=(\ebar,\overline{x})\in\Omega^{\calE,d}$, a time
$t \geq0$
and a probability measure $\Q$ on $(\Omega^{\calE,d},{\cal F}^{\sE
,d})$, let
$\Psi_{\omegabar,t}\dvtx\Omega^{\calE,d}\rightarrow\Omega^{\calE,d}$
denote the
function
%
%
\begin{equation}
\label{eq:1} \Psi_{\omegabar,t}(e,x)= \bigl(\ebar, \nabla(\overline{x},t)+x-
\nabla(x,t) \bigr),
\end{equation}
and set $\delta_{\omegabar}\otimes_t\Q=
\Q\circ\Psi^{-1}_{\omegabar,t}$. The reader can easily check that the
measure $\delta_{\omegabar}\otimes_t\Q$ is uniquely determined by the
properties
%
%
\begin{eqnarray}
\label{eq:2}
&(\delta_{\omegabar}\otimes_t\Q) \bigl[E = \ebar,
X_s = \overline{x}(s)\ \forall s \leq t\bigr] = 1,&
\\
\label{eq:3}
&(\delta_{\omegabar}\otimes_t\Q)\bigl[\Delta(X, t)
\in A\bigr] = \Q\bigl[\Delta(X, t) \in A\bigr]\qquad \forall A
\in\frakC_0^d.&
\end{eqnarray}
If $\Q[X_t=\overline{x}(t)] = 1$, $\calE= \R^d$ and we identify
$\Omega^{\calE,d}$
with $C^d$ in the natural way, then this notation reduces to the
construction given in Lemma 6.1.1 of~\cite{SV2}.

In the next step, we concatenate an initial
probability
measure and a probability kernel at a stopping time.

%
%
\begin{definition}\label{D3.4}
Let $(\Omega',\sF')$ and $(\Omega'',\sF'')$ be
measurable spaces.
We say that a function
$Q\dvtx\Omega'\times\sF''\rightarrow[0,1]$
is a {probability kernel from $(\Omega',\sF')$ to
$(\Omega'',\sF'')$} provided:
\begin{longlist}[(a)]
\item[(a)]
$Q(\omega',A'')$ is an $\sF'$-measurable
function of $\omega'\in\Omega'$
for each $A''\in\sF''$,
\item[(b)]
$Q(\omega',\cdot)$ is a probability measure on
$(\Omega'',\sF'')$ for each $\omega'\in\Omega'$.
\end{longlist}
\end{definition}

%
%
\begin{proposition}\label{P3.5}
Let $\P$ be a probability measure on
$(\Omega^{\calE,d},{\cal F}^{\sE,d})$, let $T$ be a finite (for every
$\omega$) $\{{\cal F}^{\sE,d}_t\}_{t\geq0}$-stopping time
and let $Q$ be a probability
kernel from $(\Omega^{\calE,d},{\cal F}^{\sE,d}_T)$ to
$(\Omega^{\calE,d},{\cal F}^{\sE,d})$. Then there exists a
unique probability measure $\P\otimes_TQ$ on
$(\Omega^{\calE,d},{\cal F}^{\sE,d})$ such that:
\begin{longlist}[(a)]
\item[(a)]
$\P\otimes_TQ[A]=\P[A], A\in{\cal F}^{\sE,d}_T$,
\item[(b)]
the random variable
$\omega\mapsto(\delta_{\omega}
\otimes_{T(\omega)}
Q(\omega,\cdot) )[F]$
is a version of the conditional probability
$(\P\otimes_TQ)[F | {\cal F}^{\sE,d}_T]$ for all $F\in{\cal
F}^{\sE,d}$.
\end{longlist}
\end{proposition}

\begin{pf}
When the initial condition
$Q(\omegabar,
\{E=\ebar,X_{T(\omegabar)}
=X_{T(\omegabar)}(\omegabar)\}) = 1$ holds for each
$\omegabar\in\Omega^{\calE,d}$, the result
follows in the same way as
Theorem~6.1.2 of~\cite{SV2}. To handle the
general case, we modify
the initial segment of each
path to ensure that the proper
initial condition holds.
Let $\Psi_{\omegabar,t}$ be defined as in
(\ref{eq:1}) and set
$\Qhat(\omegabar, \cdot) = Q(\omegabar,
\cdot) \circ\Psi^{-1}_{\omegabar, T(\omegabar)}$.
The map
$(\omegabar, \omega) \mapsto
\Psi_{\omegabar, T(\omegabar)}(\omega)$
is ${\cal F}^{\sE,d}_T \otimes{\cal F}^{\sE,d}/ {\cal F}^{\sE,d}$-measurable,
so $\Qhat$ is a
probability kernel from $(\Omega^{\calE,d},{\cal F}^{\sE,d}_T)$
to $(\Omega^{\calE,d},{\cal F}^{\sE,d})$.
It follows from the definition of
$\Psi$ that
$\Qhat(\omegabar,
\{E=\ebar,X^{T(\omegabar)}
=X^{T(\omegabar)}(\omegabar)\}) = 1$
for each $\omegabar\in
\Omega^{\calE,d}$, so we may apply\vspace*{1pt} the previous
case to conclude that there
exists a unique measure
$\P\otimes_T \widehat Q$ such that (a) and (b)
hold when $Q$ is replaced with $\widehat Q$. But
the operator $\omega
\mapsto\Psi_{\omegabar,t}(\omega)$ is idempotent, so
$\delta_{\omegabar} \otimes_{T(\omegabar)}
Q(\omegabar,\cdot) =
\delta_{\omegabar} \otimes_{T(\omegabar)}
\widehat Q(\omegabar,\cdot)$
for each $\omegabar\in\Omega^{\calE,d}$,
and $\P\otimes_T Q = \P\otimes_T
\widehat Q$ is in fact the unique
measure which satisfies (a) and (b).
\end{pf}

We now begin concatenating probability measures.

%
%
\begin{corollary}[(Two-fold concatenation)]\label{C3.6}
Let $\P_1$ and $\P_2$ be probability measures
on $\Omega^{\calE,d}$, let $T$ be a finite (for every
$\omega$) $\{{\cal F}^{\sE,d}_t\}_{t\geq0}$-stopping time,
let $\sG$ be a sub-$\sigma$-field of ${\cal F}^{\sE,d}_T$
and assume that $\P_1|{\sG}\ll\P_2|_{\sG}$.
Then there exists a unique measure, denoted
$\P_1\otimes_{T,\sG}\P_2$, such that:
\begin{longlist}[(a)]
\item[(a)]
$\P_1\otimes_{T,\sG}\P_2[A]=\P_1[A],
A\in{\cal F}^{\sE,d}_T$,
\item[(b)]
for every set $B\in\sG\vee\sigma(\Delta(X,T))$,
every version of $\P_2[B | \sG]$ is a
version of
$(\P_1\otimes_{T,\sG}\P_2)[B | {\cal F}^{\sE,d}_T]$,
\item[(c)]
if $\P_1$ and $\P_2$ agree on $\sG$, then
$\P_1\otimes_{T,\sG}\P_2$ and $\P_2$ agree on
$\sG\vee\sigma(\Delta(X,T))$.
\end{longlist}
\end{corollary}

\begin{pf}
Because $\Omega^{\calE,d}$ is a Polish space, there
exists a $\sG$-measurable probability kernel
$Q$ from
$(\Omega^{\calE,d},{\cal F}^{\sE,d})$ to
$(\Omega^{\calE,d},{\cal F}^{\sE,d})$ such that for every
$F\in{\cal F}^{\sE,d}$, $Q(\cdot,F)$ is a
version of $\P_2[F | \sG]$
(\cite{SV2}, Theorem 1.1.6). Using
Proposition~\ref{P3.5}, we define
$\P_1\otimes_{T,\sG}\P_2\define\P_1\otimes_TQ$.
Property (a) of the
corollary is property (a) of Proposition~\ref{P3.5}.

Given $\omega\in\Omega^{\calE,d}$ and $F \in{\cal F}^{\sE,d}$,
set $\Qhat(\omega, F) = (\delta_\omega
\otimes_{T(\omega)}Q(\omega,\cdot))[F]$.
Property (b) of Proposition
\ref{P3.5} asserts that $\Qhat(\cdot, F)$
is a version of
$(\P_1\otimes_TQ)[F | {\cal F}^{\sE,d}_T]$ for all
$F \in{\cal F}^{\sE,d}$.
Galmarino's test (\cite{DM},
Theorem IV.100) for the
filtered space $(\Omega^{\calE,d},{\cal F}^{\sE,d},
\{{\cal F}^{\sE,d}_t\}_{t\geq0})$
says that $E(\omega)=E(\omegabar)$
and $X_u(\omega)=X_u(\omegabar)$ for
$0\leq u\leq T(\omegabar)$ imply
$Y(\omega)=Y(\omegabar)$ for every
$\sF_T^{\sE,d}$-measurable random variable~$Y$.
In particular, if $A\in\sF_T^{\sE,d}$, then
$\omega\in A$ if and only if $ (E(\omega),
X^{T(\omega)}(\omega) ) \in A$.
Therefore, $\Qhat(\omega, A \cap F) = \ind_A(\omega)
\Qhat(\omega, F)$ for all $\omega\in\Omega^{\calE,d}$,
$A \in{\cal F}^{\sE,d}_T$ and
$F \in{\cal F}^{\sE,d}$ by (\ref{eq:2}).
If $B=A\cap\{\Delta(X,T)\in D\}$
with $A\in\sG$ and $D\in\frakC_0^d$,
then (\ref{eq:3}) implies
%
%
\begin{equation}
\label{eq:4} \Qhat(\cdot, B) = \ind_A Q\bigl(\cdot, \bigl\{
\Delta(X, T) \in D\bigr\}\bigr) = \P_2[B | \sG],\qquad \P_2
\mbox{-a.s.}
\end{equation}
It then follows from Dynkin's $\pi$--$\lambda$
theorem that (\ref{eq:4}) holds for
all $B \in\sG\vee\sigma(\Delta(X,T))$.
From (a) we have
$\P_1\otimes_{T,\sG}\P_2|_{\sG}=\P_1|_{\sG}$,
and we have assumed
$\P_1|_{\sG}\ll\P_2|_{\sG}$,
so the fact that $\Qhat(\cdot, B)$ is a version
of both $(\P_1\otimes_{T,\sG}\P_2)[B | {\cal F}^{\sE,d}_T]$
and $\P_2[B | \sG]$
implies that every version of
$\P_2[B | \sG]$ is also a version of
$(\P_1\otimes_{T,\sG}\P_2)[B | {\cal F}^{\sE,d}_T]$.

For (c), assume that $\P_1$ and $\P_2$ agree
on $\sG$. Property (a) implies that
$\P_1|_{\sG}=(\P_1\otimes_{T,\sG}\P_2)|_{\sG}$,
and hence,
$\P_2|_{\sG}=\P_1\otimes_{T,\sG}\P_2|_{\sG}$.
For $B\in\sG\vee\sigma(\Delta(X,T))$,
we have from (b) that
$\P_2[B | \sG]
=(\P_1\otimes_{T,\sG}\P_2)[B | \sG]$,
and we can integrate both sides over $\Omega^{\calE,d}$
with respect to
$\P_2|_{\sG}=\P_1\otimes_{T,\sG}\P_2|_{\sG}$
to obtain (c).

Uniqueness of $\P_1\otimes_{T,\sG}\P_2$
follows from the fact that (b) specifies this
measure on $\sG\vee\sigma(\Delta(X,T))$
conditioned on ${\cal F}^{\sE,d}_T$, up to
$\P_2|_{\sG}$-equivalence. Furthermore, (a)
specifies this measure to be $\P_1$ on ${\cal F}^{\sE,d}_T$,
and hence, on $\sG$.
But $\P_1\ll\P_2$, and hence, the integral
in the equation
$(\P_1\otimes_{T,\sG}\P_2)[B]
=\int_{\Omega^{\calE,d}}\P_2[B | \sG] \,d\P_1$ for
$B\in\sG\vee\sigma(\Delta(X,T))$
is well defined.
We see then that properties (a) and (b)
specify the measure $\P_1\otimes_{T,\sG}\P_2$ on
$\sG\vee\sigma(\Delta(X,T))$
and on ${\cal F}^{\sE,d}_T$. These
two $\sigma$-fields generate ${\cal F}^{\sE,d}$, and thus
the measure is uniquely determined
on ${\cal F}^{\sE,d}$ by properties (a) and~(b).
\end{pf}

%
%
\begin{proposition}[(Three-fold concatenation)]\label{P3.7}
Let\vspace*{1pt} $\P_1$, $\P_2$ and $\P_3$ be probability
measures on $\Omega^{\calE,d}$ and let $0\leq S\leq T$ be
finite (for every $\omega$)\break
$\{{\cal F}^{\sE,d}_t\}_{t\geq0}$-stopping times. Let
$\sG$ be a sub-$\sigma$-field of ${\cal F}^{\sE,d}_S$
and let $\sH$ be a sub-$\sigma$-field of
$\sG\vee\sigma(\Delta(X^T,S))$, which is a
sub-$\sigma$-field of ${\cal F}^{\sE,d}_T$. Assume that
$T-S$ is
$\sG\vee\sigma(\Delta(X,S))$-measurable.
If $\P_1|_{\sG}\ll\P_2|_{\sG}$ and
$\P_2|_{\sH}\ll\P_3|_{\sH}$, then:
\begin{longlist}[(a)]
\item[(a)]
$\P_1|_{\sG}\ll
(\P_2\otimes_{T,\sH}\P_3)|_{\sG}$,
\item[(b)]
$(\P_1\otimes_{S,\sG}\P_2)|_{\sH}
\ll\P_3|_{\sH}$,
\end{longlist}
so that both
$\P_1\otimes_{S,\sG}(\P_2\otimes_{T,\sH}\P_3)$
and
$(\P_1\otimes_{S,\sG}\P_2)\otimes_{T,\sH}\P_3$
are defined, and
\begin{longlist}[(c)]
\item[(c)]
$\P_1\otimes_{S,\sG}(\P_2\otimes_{T,\sH}\P_3)
=(\P_1\otimes_{S,\sG}\P_2)\otimes_{T,\sH}\P_3$.
\end{longlist}
\end{proposition}

\begin{pf}
We simplify notation by writing
$\P_{12}=\P_1\otimes_{S,\sG}\P_2$,
$\P_{23}=\P_2\otimes_{T,\sH}\P_3$,
$\P_{1,23}=\P_1\otimes_{S,\sG}
(\P_2\otimes_{T,\sH}\P_3)$ and
$\P_{12,3}=(\P_1\otimes_{S,\sG}\P_2)
\otimes_{T,\sH}\P_3$.
For (a), we note from Corollary~\ref{C3.6}(a)
that $\P_{23}$ agrees with $\P_2$ on ${\cal F}^{\sE,d}_T$,
and hence, on $\sG$. Property (a) follows from
$\P_1|_{\sG}\ll\P_2|_{\sG}$.
For (b), let $A\in\sH$ satisfy $\P_3[A]=0$. By
assumption, we also have $\P_2[A]=0$, and hence,
$0$ is a version of $\P_2[A | \sG]$. Being in
$\sH$, $A$ is also in
$\sG\vee\sigma(\Delta(X,S))$, and according to
Corollary~\ref{C3.6}(b), $0$ is a version of
$\P_{12}[A | {\cal F}^{\sE,d}_S]$. Therefore, $\P_{12}[A]=0$.

The collection of sets of the form
$A\cap B\cap C$, where
$A\in{\cal F}^{\sE,d}_S$, $B\in\break\sigma(\Delta(X^T,S))$
and $C\in\sigma(\Delta(X,T))$, is closed under
finite intersections and generates ${\cal F}^{\sE,d}$.
Thus, to prove (c), it suffices to show
that the desired equation holds when both sides
are evaluated for a set of this form. Let $A$,
$B$ and $C$ be as described, and let $G$ be
in $\sG$. Let $Z$ be a version of
$\E_3[\ind_C | \sH]$ and $Y$ a version of
$\E_2[\ind_BZ | \sG]$. Corollary~\ref{C3.6}(b)
implies that $Z$ is a version of
$\P_{23}[C | {\cal F}^{\sE,d}_T]$. This, combined with
Corollary~\ref{C3.6}(a), implies
\[
\E_{23}[\ind_GY] =\E_2[\ind_GY]
=\E_2[\ind_{G\cap B}Z] =\E_{23}[
\ind_{G\cap B}Z] =\E_{23}[G\cap B\cap C].
\]
We see then that
$Y=\E_2[\ind_B\E_3[\ind_C | \sH] | \sG]$
is a version of $\E_{23}[\ind_{B\cap C} | \sG]$,
a fact we use along with repeated applications
of Corollary~\ref{C3.6}(a), (b) and (c) in the
chain of equalities
\begin{eqnarray*}
\P_{1,23}[A\cap B\cap C] &=& \E_{1,23} \bigl[
\ind_A\E_{1,23}\bigl[\ind_{B\cap C} | {\cal
F}^{\sE,d}_S\bigr] \bigr]
\\
&=& \E_{1,23} \bigl[\ind_A\E_{23}[
\ind_{B\cap C} | \sG] \bigr]
\\
&=& \E_{1,23} \bigl[\ind_A\E_2 \bigl[
\ind_B\E_3[\ind_C | \sH] | \sG\bigr] \bigr]
\\
&=& \E_1 \bigl[\ind_A\E_2 \bigl[
\ind_B\E_3[\ind_C | \sH] | \sG\bigr] \bigr]
\\
&=& \E_{12} \bigl[\ind_A\E_{12} \bigl[
\ind_B\E_3 [\ind_C | \sH] | {\cal
F}^{\sE,d}_S \bigr] \bigr]
\\
&=& \E_{12} \bigl[\ind_{A\cap B} \E_3[
\ind_C | \sH] \bigr]
\\
&=& \E_{12,3} \bigl[\ind_{A\cap B}\E_{12,3} \bigl[
\ind_C | {\cal F}^{\sE,d}_S\bigr] \bigr]
\\
&=& \P_{12,3}[A\cap B\cap C].
\end{eqnarray*}
\upqed
\end{pf}

\begin{pf*}{Proof of Theorem~\ref{T3.3}}
Let $m$ satisfy $0\leq m\leq n-1$.
According to Definition~\ref{D3.2},
\[
\sG_{m+1}\subset\sH_{m+1} =\sG_m\vee\sigma
\bigl(\Delta\bigl(X^{T_{m+1}},T_m\bigr) \bigr) \subset
\sG_m\vee\sigma\bigl(\Delta(X,T_m) \bigr).
\]
If $0\leq m\leq n-2$, we further have
\[
\sG_{m+2}\subset\sH_{m+2}=\sG_{m+1}\vee\sigma
\bigl(\Delta\bigl(X^{T_{m+2}},T_{m+1}\bigr) \bigr) \subset
\sG_m\vee\sigma\bigl(\Delta(X,T_m) \bigr).
\]
Iterating this process, we obtain the relation
$\sG_{j}\subset\sG_m\vee
\sigma(\Delta(X,T_m) )$ for
$j=m,m+1,\ldots,n$.
Consequently,
%
%
\begin{equation}
\label{3.10} \sG_j\vee\sigma\bigl(\Delta\bigl(X^{T_{j+1}},T_j
\bigr) \bigr) \subset\sG_m\vee\sigma\bigl(\Delta(X,T_m)
\bigr),\qquad 0\leq m\leq j\leq n.
\end{equation}

We now proceed by induction on $m$. The induction hypothesis
corresponding to~$m$, where $m=0,\ldots,n$, is the existence of a
measure $\P^m$ such that:
\begin{longlist}[(ii)]
\item[(i)]
$\P^m[A]=\P[A]$ for all $A\in\sH_i$ and
$0\leq i\leq n+1$,
\item[(ii)]
for $B\in\sH_{i+1}$ and $0\leq i\leq m-1$,
every $\P$-version of $\P[B | \sG_i]$ is a
$\P^m$-version of $\P^m[B | {\cal F}^{\sE,d}_{T_i}]$.
\end{longlist}
The base case is $\P^0=\P$, a case for which (i)
trivially holds and (ii) is vacuous.

Assume the induction hypothesis for some integer $m$. Because
$\sG_{m}\subset\sH_{m}$ and the measures $\P^m$ and $\P$ agree on
$\sH_{m}$, we may invoke Corollary~\ref{C3.6} to define
$\P^{m+1}\define\P^m\otimes_{T_{m},\sG_{m}}\P$. If $A\in\sH_j$ for some
$j$, $0\leq j\leq m$, then $A\in{\cal F}^{\sE,d}_{T_{m}}$ and
$\P^{m+1}[A]=\P^{m}[A]=\P[A]$ by Corollary~\ref{C3.6}(a) and part (i)
of the induction hypothesis. If $m\leq j\leq n$, then (\ref{3.10})
implies
%
%
\begin{equation}
\label{HinG} \sH_{j+1}\define\sG_j\vee\sigma\bigl(\Delta
\bigl(X^{T_{j+1}},T_j\bigr) \bigr) \subset\sG_{m}\vee
\sigma\bigl(\Delta(X,T_{m}) \bigr).
\end{equation}
But Corollary~\ref{C3.6}(c) implies that $\P^{m+1}$ agrees with $\P$ on
$\sG_m\vee\sigma(\Delta(X,T_m) )$. Hence, $\P^{m+1}$ satisfies (i).

For some $i$, $0\leq i\leq m-1$, let $B\in\sH_{i+1}$ and $A\in{\cal
F}^{\sE,d}_{T_i}$ be given. Suppose $Z$ is a version of $\P[B |
\sG_i]$, so that both $A\cap B$ and $\ind_AZ$ are ${\cal
F}^{\sE,d}_{T_m}$-measurable\vspace*{1pt} [recall (\ref{3.1a})]. Corollary
\ref{C3.6}(a) (used twice) and part (ii) of the induction hypothesis
imply
\[
\E^{m+1}[\ind_AZ]=\E^m[\ind_AZ]=
\P^m[A\cap B] =\P^{m+1}[A\cap B],
\]
showing that $Z$ is a version of
$\P^{m+1}[B | {\cal F}^{\sE,d}_{T_i}]$. Finally, suppose
$B$ is in $\sH_{m+1}$, which is a
sub-$\sigma$-field of
$\sG_m\vee\sigma(\Delta(X,T_m))$.
Corollary~\ref{C3.6}(b) says that every version
of $\P[B | \sG_m]$ is a version of
$\P^{m+1}[B | {\cal F}^{\sE,d}_{T_m}]$. This establishes
(ii) with $m+1$ replacing $m$.

The induction argument above constructs
$\P^{\otimes\Pi}\define\P^{n+1}$ that satisfies
(\ref{3.2}) and (\ref{3.3}). To see that this
measure is unique, we show that (\ref{3.2}) and
(\ref{3.3}) determine its value on sets of the
form $\bigcap_{i=0}^{n+1}B_i$, where
$B_0\in{\cal F}^{\sE,d}_0=\sH_0$
and
$B_{i+1}\in\sigma(\Delta(X^{T_{i+1}},T_i))
\subset\sH_{i+1}$
for $i=0,\ldots,n$. This collection of sets
is closed under finite intersections and
generates ${\cal F}^{\sE,d}$. For such a set, repeated
application of (\ref{3.3}), followed
by a final application of (\ref{3.2}), yields
%
%
\begin{eqnarray}\label{3.11}\quad
&&
\P^{\otimes\Pi} \Biggl[\bigcap_{i=0}^{n+1}B_i
\Biggr]
\nonumber
\\
&&\qquad= \E^{n+1} \bigl[\ind_{B_0}\E^{n+1} \bigl[
\ind_{B_1}\cdots\E^{n+1} \bigl[ \ind_{B_n}
\E^{n+1}\bigl[\ind_{B_{n+1}} | {\cal F}^{\sE,d}_{T_n}
\bigr] | {\cal F}^{\sE,d}_{T_{n-1}} \bigr] \cdots| {\cal
F}^{\sE,d}_0 \bigr] \bigr]
\nonumber
\\
&&\qquad= \E^{n+1} \bigl[\ind_{B_0}\E^{n+1} \bigl[
\ind_{B_1}\cdots\E^{n+1} \bigl[ \ind_{B_n}\E[
\ind_{B_{n+1}} | \sG_n] | {\cal F}^{\sE,d}_{T_{n-1}}
\bigr] \cdots| {\cal F}^{\sE,d}_0 \bigr] \bigr]
\nonumber
\\
&&\qquad= \E^{n+1} \bigl[\ind_{B_0}\E^{n+1} \bigl[
\ind_{B_1}\cdots\E\bigl[\ind_{B_n}\E[\ind_{B_{n+1}} |
\sG_n] \sG_{n-1} \bigr] \cdots| {\cal
F}^{\sE,d}_0 \bigr] \bigr]
\\
&&\hspace*{2.6pt}\qquad\vdots
\nonumber
\\
&&\qquad= \E^{n+1} \bigl[\ind_{B_0}\E\bigl[\ind_{B_1}
\cdots\E\bigl[\ind_{B_n}\E[\ind_{B_{n+1}} | \sG_n]
| \sG_{n-1} \bigr] \cdots| \sG_0 \bigr] \bigr]
\nonumber
\\
&&\qquad= \E\bigl[\ind_{B_0}\E\bigl[\ind_{B_1}\cdots\E\bigl[
\ind_{B_n}\E[\ind_{B_{n+1}} | \sG_n] |
\sG_{n-1} \bigr] \cdots| \sG_0 \bigr] \bigr].\nonumber
\end{eqnarray}
The proof of Theorem
\ref{T3.3} is complete.
\end{pf*}

%
%
\begin{remark}\label{R3.8}
We see from the proof of Theorem~\ref{T3.3} that
\[
\P^{\otimes\Pi} =\P\otimes_{T_0,\sG_0}\P\otimes_{T_1,\sG_1} \cdots
\otimes_{T_n,\sG_n}\P,
\]
where the associative property of Proposition
\ref{P3.7}(c) makes the grouping of the
$\otimes_{T_i,\sG_i}$ operators irrelevant.
Equation (\ref{3.11}) provides insight into
the nature of $\P^{\otimes\Pi}$. If $\sG_i$
is equal to ${\cal F}^{\sE,d}_{T_i}$ for each $i$, then the
last iterated conditional expectation in
(\ref{3.11}) collapses to
$\P[\bigcap_{i=0}^{n+1}B_i]$, and
$\P^{\otimes\Pi}$ agrees with $\P$. At the
other extreme, if $\sG_i$ is the trivial
$\sigma$-field $\{\varnothing,\Omega^{\calE,d}\}$ for each
$i$, then this iterated conditional expectation
becomes $\prod_{i=0}^{n+1}\P[B_i]$, and
increments of the path fragments over
$[T_i,T_{i+1}]$ are independent of one another
under $P^{\otimes\Pi}$ but have the same unconditional
distribution as under $\P$.
\end{remark}

\subsection{Properties preserved by
concatenation}\label{SPropPres}

%
%
\begin{proposition}\label{P3.9}
Let $\P$ be a probability measure and let
$(T_i,\sG_i)_{i=0}^n$ be an extended
partition on
$(\Omega^{\calE,d},{\cal F}^{\sE,d})$. Let $A$ be an $\{{\cal
F}^{\sE,d}_t\}_{t\geq0}$-adapted continuous
real-valued process on $\Omega^{\calE,d}$, and assume
that $\Delta(A,T_i)$ is
$\sG_i\vee\sigma(\Delta(X,T_i))$-measurable for
$i=0,\ldots,n$.
\begin{longlist}[(a)]
\item[(a)]
The total variation of $A$ on $[0,\infty)$
is $\P$-almost surely finite if and only if
it is $\P^{\otimes\Pi}$-almost surely finite.
\item[(b)]
The process $A$ is $\P$-almost surely absolutely
continuous if and only if it
is $\P^{\otimes\Pi}$-almost surely absolutely continuous.
\end{longlist}
\end{proposition}

\begin{pf}
We set $\P^0=\P$ and
$\P^{i+1}=\P^i\otimes_{T_i,\sG_i}\P$,
$i=1,\ldots,n$.
Then $\P^{\otimes\Pi}=\P^{n+1}$.
For (a), we proceed by induction
on $i=0,1,\ldots,n$, assuming that:
\begin{longlist}[(a$_i$)]
\item[(a$_i$)]
the total variation of $A$ on $[0,\infty)$
is $\P$-almost surely finite if and only if
it is $\P^i$-almost surely finite.
\end{longlist}
On ${\cal F}^{\sE,d}_{T_i}$, the probability measures $\P^i$
and $\P^{i+1}$ agree [Corollary~\ref{C3.6}(a)],
and hence, $A$ restricted
to $[0,T_i]$ is $\P^i$-a.s. of finite
variation if and only if $A$ restricted to
$[0,T_i]$ is $\P^{i+1}$-a.s. of finite
variation. The variation of $A$
on subintervals in $[T_i,\infty)$ is a function
of $\Delta(A,T_i)$, which is
$\sG_i\vee\sigma(\Delta(X,T_i))$-measurable,
and on this $\sigma$-field, the measures $\P^i$
and $\P^{i+1}$ agree [Corollary~\ref{C3.6}(c)].
Therefore, $A$ restricted to $[T_i,\infty)$
is $\P^i$-a.s. of finite variation if and only
if $A$ restricted to $[T_i,\infty)$ is
$\P^{i+1}$-a.s. of finite variation. We
conclude that
$A$ has finite total variation on $[0,\infty)$
$\P^i$-almost surely if and only if it has
finite total variation $\P^{i+1}$-almost surely.
Combining this with the
induction hypothesis (a$_i$), we obtain
the induction hypothesis with $i+1$ replacing~$i$.

The continuous process $A$ is absolutely
continuous on $[0,\infty)$ if and only if it is
absolutely continuous on $[0,T_i]$ and
absolutely continuous on $[T_i,\infty)$.
Therefore, we can imitate the proof of (a)
to obtain (b).
\end{pf}

%
%
\begin{proposition}\label{P3.10}
Let $\P$ be a probability measure and let
$(T_i,\sG_i)_{i=0}^n$ be an extended
partition on
$(\Omega^{\calE,d},{\cal F}^{\sE,d})$. Let $A$ be an
$\{{\cal F}^{\sE,d}_t\}_{t\geq0}$-adapted continuous
$\R^d$-valued process on $\Omega^{\calE,d}$
with $A_0=0$,
and assume that $\Delta(A,T_i)$ is
$\sG_i\vee\sigma(\Delta(X,T_i))$-measurable for
$i=0,\ldots,n$.
Assume there exists a measurable $\R^d$-valued
process $\alpha$ such that the set
%
%
\begin{equation}\quad
\label{3.14} J(\omega)\define\biggl\{t\in[0,\infty)\dvtx\frac{\partial
}{\partial t}A_t(
\omega) \mbox{ exists but is not equal to } \alpha_t(\omega) \biggr\}
\end{equation}
has Lebesgue measure zero for
$\P$-almost every and $\P^{\otimes\Pi}$-almost every
$\omega\in\Omega^{\calE,d}$. Then
%
%
\begin{equation}
\label{3.15} \P\biggl[A_t=\int_0^t
\alpha_u\,du\ \forall t\in[0,\infty) \biggr]=1
\end{equation}
if and only if
%
%
\begin{equation}
\label{3.16} \P^{\otimes\Pi} \biggl[A_t=\int_0^t
\alpha_u\,du\ \forall t\in[0,\infty) \biggr]=1.
\end{equation}
When the equalities (\ref{3.15}) and
(\ref{3.16}) hold, we also have
%
%
\begin{equation}
\label{3.17} \E\int_0^{S}f(\alpha_u)
\,du =\E^{\otimes\Pi}\int_0^{S} f(
\alpha_u) \,du
\end{equation}
for every nonnegative, Borel-measurable function
$f\dvtx\R^d\rightarrow\R$ and
$\{{\cal F}^{\sE,d}_t\}_{t\geq0}$-stopping
time $S$ satisfying
$(S-T_i)^+\in\sG_i\vee\sigma(\Delta(X,T_i))$
for $i=0,\ldots,n$.
\end{proposition}

\begin{pf}
Assume (\ref{3.15}).
Then each component of $A$ is $\P$-a.s. absolutely continuous.
Proposition~\ref{P3.9} implies that the components of $A$ are
$\P^{\otimes\Pi}$-a.s. absolutely continuous
as well. Therefore, for
$\P^{\otimes\Pi}$-almost every $\omega$, the set
\[
C(\omega)\define\biggl\{t\in[0,\infty)\dvtx\frac{\partial}{\partial t}A_t(
\omega) \mbox{ exists} \biggr\}
\]
has full Lebesgue measure, and by the assumption
about $J(\omega)$, the set
\[
D(\omega)\define\biggl\{t\in[0,\infty)\dvtx\frac{\partial}{\partial t}A_t(
\omega) \mbox{ exists and is equal to } \alpha_t(\omega) \biggr\}
\]
also has full Lebesgue measure for
$\P^{\otimes\Pi}$-almost
every $\omega$. This implies (\ref{3.16}).
This argument is reversible; (\ref{3.16})
implies (\ref{3.15}).

We now assume (\ref{3.15}) and (\ref{3.16}). The
$\sG_i\vee\sigma(\Delta(X,T_i))$-measurability of $\Delta(A,T_i)$
together with the ${\cal F}^{\sE,d}_{T_{i+1}}$-measurability of
$\Delta(A^{T_{i+1}},T_i)$ implies the $\sH_{i+1}$-measurability of
$\Delta(A^{T_{i+1}},T_i)$. Because $A$ is adapted and continuous,
$A_{T_i(\omega)+t}(\omega) \ind_{\{0\leq t<
T_{i+1}(\omega)-T_i(\omega)\}}$ is a jointly $\sH_{i+1}\otimes{\cal
B}[0,\infty)$-measurable function of $(\omega,t)$, where
$\sB[0,\infty)$ is the Borel $\sigma$-field on $[0,\infty)$ (recall
Remark~\ref{R3.2}). The same is then true for the right-hand derivative
$\frac{\partial^+}{\partial t}A_{T_i(\omega)+t} \ind_{\{0\leq t<
T_{i+1}(\omega) -T_i(\omega)\}}$, where we set this right-hand
derivative equal to an arbitrary value whenever the limit of the
relevant difference quotient does not exist. By assumption, $(S-T_i)^+$
is also $\sG_i\vee\sigma(\Delta(X,T_i))$-measurable. Therefore,
$(T_{i+1}-T_i)\wedge(S-T_i)^+ =T_{i+1}\wedge S-T_i\wedge S$ is
$\sH_{i+1}$-measurable. But on each $\sH_{i+1}$, the measures $\P$ and
$\P^{\otimes\Pi}$ agree, which implies that for every nonnegative
Borel-measurable function $f\dvtx\R^d\rightarrow\R$,
\[
\E\int_0^{T_{i+1}\wedge S-T_i\wedge S} f(\alpha_{T_i+u}) \,du =
\E^{\otimes\Pi} \int_0^{T_{i+1}\wedge S-T_i\wedge S} f(
\alpha_{T_i+u}) \,du.
\]
Summing over $i=0,1,\ldots,n$,
we obtain (\ref{3.17}).
\end{pf}

%
%
\begin{example}[(Example~\ref{SimpExPart1} continued)]
Consider the extended partition and
probability measures $\P$ and $\P^{\otimes\Pi}$
of Example~\ref{SimpExPart1}.
We take $A=X$ so that $A_0=0$ and
$\Delta(A,T_i)$ is
$\sG_1\vee\sigma(\Delta(X,T_i))$-measurable
for $i=0,1$.
We define the adapted processes
\begin{eqnarray*}
\alpha_t(\omega) &=& \biggl\{\limsup_{\varepsilon\downarrow0} 0 \vee
\frac{\omega(\varepsilon)}{\varepsilon} \wedge1 \biggr\} \ind_{(0,\infty]}(t),
\\
\beta_t(\omega) &=& \biggl\{\limsup_{\varepsilon\downarrow0} 0 \vee
\frac{\omega(\varepsilon)}{\varepsilon} \wedge1 \biggr\} \ind
_{(0,1]}(t) \\
&&{}+ \biggl\{
\limsup_{\varepsilon\downarrow0} 0 \vee\frac{\omega(1+\varepsilon
)-\omega(1)}{\varepsilon} \wedge1 \biggr\}
\ind_{(1,\infty]}(t)
\end{eqnarray*}
and the sets $E = \{A_t=\int_0^t\alpha_u \,du\ \forall
t\in[0,\infty)\}$ and $F = \{A_t=\int_0^t\beta_u \,du\ \forall
t\in[0,\infty)\}$. Then
we have $\P[E] = \P[F] = \P^{\otimes\Pi}[F] = 1$, but
$\P^{\otimes\Pi}[E] = 1/2$.

If we let $K(\omega)$ denote the set obtained by replacing $\alpha$
with $\beta$ in (\ref{3.14}), then we see that $K(\omega)$ is a
Lebesgue-null set $\P$-almost surely and $\P^{\otimes\Pi}$-almost
surely. On the other hand,
$J(\omega)$ defined by (\ref{3.14})
is a Lebesgue-null set $\P$-almost surely, but has
strictly positive Lebesgue measure with
strictly positive
$\P^{\otimes\Pi}$-probability. In particular,
we see that
(\ref{3.15}) and (\ref{3.16}) may
not be equivalent in this situation.
\end{example}

%
%
\begin{corollary}\label{C3.11}
Let $\P$ be a probability measure on
$(\Omega^{\calE,d},{\cal F}^{\sE,d})$ and for each positive integer
$m$, let
$\Pi^m\define(T_i^m,\sG_i^m)_{i=0}^{N(m)}$
be an extended partition. Let $A$ be an
$\{{\cal F}^{\sE,d}_t\}_{t\geq0}$-adapted
continuous $\R^d$-valued process on
$\Omega^{\calE,d}$ with $A_0=0$, and assume that
$T^m_i$ and $\Delta(A,T_i^m)$ are
$\sG_i^m\vee\sigma(\Delta(X,T_i^m))$-measurable
for $i=1,\ldots,N(m)$ and $m=1,2,\ldots\,$.
Let $\alpha$ be a measurable $\R^d$-valued
process such that
$A_t=\int_0^t\alpha_u\,du$ for
\textit{every} $t\geq0$, $\P$-almost surely,
and assume that the set $J(\omega)$
defined by (\ref{3.14}) has Lebesgue
measure zero for \textit{every} $\omega\in\Omega^{\calE,d}$.
Finally, assume
%
%
\begin{equation}
\label{3.18b} \E\int_0^t\|\alpha_u
\| \,du<\infty,\qquad t\geq0.
\end{equation}
Then the following hold.
\begin{longlist}[(a)]
\item[(a)]
For every $t\in[0,\infty)$,
$\alpha$ restricted to $[0,t]$ is uniformly
integrable with respect to
the collection of product measures
$\{\P^{\otimes\Pi^m}\times
\lambda_{[0,t]}\}_{m=1}^{\infty}$,
where $\lambda_{[0,t]}$ denotes Lebesgue
measure on $[0,t]$.
\item[(b)]
The collection of measures
$\{\P^{\otimes\Pi^m}\circ
A^{-1}\}_{m=1}^{\infty}$
on $C_0^d$ is tight.
\end{longlist}
\end{corollary}

\begin{pf}
For (a), fix $t\in[0,\infty)$. Given $\varepsilon>0$, (\ref{3.18b})
guarantees that there exists $M_\varepsilon>0$ so large that
$\E\int_0^t\|\alpha_u\| \ind_{\{\|\alpha_u\|\geq M_{\varepsilon}\}}
\,du \leq\varepsilon$. Applying Proposition~\ref{P3.10} with
$f(x)=\|x\|\ind_{\{\|x\|\geq M_{\varepsilon}\}}$ and $S=t$, we
obtain\break
$\E^{\otimes\Pi^m} \int_0^t\|\alpha_u\|\ind_{\{\|\alpha_u\|\geq
M_{\varepsilon}\}} \,du \leq\varepsilon$ for all~$m$.

For (b) it suffices to verify that for every $\varepsilon>0$, there
exists a set $\Omega_{\varepsilon}\in{\cal F}^{\sE,d}$ such that
$\P^{\otimes\Pi^m}(\Omega_{\varepsilon})\geq 1-\varepsilon$ for every
$m$ and
%
%
\begin{equation}
\label{3.19b} \lim_{\delta\downarrow0} \sup_{\omega\in\Omega
_{\varepsilon}} \sup_{0\leq s\leq v\leq t\wedge(s+\delta)}
\bigl\|A_v(\omega)-A_s(\omega)\bigr\|=0,\qquad t\geq0.
\end{equation}
Fix $\varepsilon>0$ and let $\{t_n\}_{n=1}^{\infty}$ be an
increasing sequence of positive numbers with
$\lim_{n\rightarrow\infty}t_n=\infty$.
For fixed $n$, we construct $\Omega_n$ such that
$\P^{\otimes\Pi^m}(\Omega_n)\geq1-2^{-n}\varepsilon$ for every $m$ and
%
%
\begin{equation}
\label{3.20b} \lim_{\delta\downarrow0}\sup_{\omega\in\Omega_n} \sup
_{0\leq s\leq v\leq t_n\wedge(s+\delta)}
\bigl\|A_v(\omega)-A_s(\omega)\bigr\|=0.
\end{equation}
Then $\Omega_{\varepsilon}=\bigcap_{n=1}^{\infty}\Omega_n$
satisfies (\ref{3.19b}) and $\P^{\otimes\Pi^m}(\Omega_{\varepsilon
})\geq1-\varepsilon$
for every $m$.

We fix $n$ and construct $\Omega_n$ by working through the proof
of the Borel--Cantelli lemma.
For each positive integer $k$,
part (a) implies the existence of $\delta_k>0$
for which
\[
\E^{\otimes\Pi^m} \Bigl[ \sup_{0\leq s\leq v\leq t_n\wedge(s+\delta_k)}
\|A_v-A_s
\| \Bigr]\leq\E^{\otimes\Pi^m} \biggl[\sup_{0\leq s\leq v\leq t_n\wedge
(s+\delta_k)} \int
_s^v\|\alpha_u\| \,du \biggr] \leq
2^{-2k}
\]
for all $m$. We define\vspace*{1pt} $F_k= \{\sup_{0\leq s\leq v\leq
t_n\wedge(s+\delta_k)} \int_s^v\|\alpha_u\| \,du\geq2^{-k} \}$, and
note from Chebyshev's inequality that
$\P^{\otimes\Pi^m}(F_k)\leq2^{-k}$ for every $m$ and $k$. Choose $j$
such that $2^{-(j-1)}\leq2^{-n}\varepsilon$ and set
$\Omega_n=\bigcap_{k\geq j}F_k^c$. We have
$\P^{\otimes\Pi^m}(\Omega_n^c)\leq\sum_{k=j}^{\infty}\P^{\otimes\Pi^m}(F_k)
\leq2^{-n}\varepsilon$ for every $m$, as desired. Also,
$\omega\in\Omega_n$ implies that $\sup_{0\leq s\leq v\leq
t_n\wedge(s+\delta_k)} \int_s^v\|\alpha_u\| \,du\leq2^{-k}$ for all
$k\geq j$, and hence, (\ref{3.20b}) holds.
\end{pf}

%
%
\begin{definition}\label{D3.11}
Let $Y$ be an adapted continuous $\R^d$-valued process
defined on a filtered probability
space $(\Omega,\sF,\{\sF_t\}_{t\geq0},\P)$, let $B$
be an adapted continuous $\R^d$-valued process whose components
are of finite variation and for which $B_0=0$ and let
$C$ be an adapted continuous,
$d\times d$-matrix-valued process whose components
are of finite variation and for which $C_0=0$.
We further assume that outside a $\P$-null set
that does not depend on $s$ and $t$,
the increment $C_t-C_s$ is positive semidefinite whenever
$0\leq s<t<\infty$.
We say that \textit{$Y$ is a semimartingale
with characteristic pair $(B,C)$} if the components of
$Y-B$ and $(Y-B)(Y-B)^{tr}-C$ are local martingales.
\end{definition}

%
%
\begin{proposition}\label{P3.12}
Let $\P$ be a probability measure and let
$(T_i,\sG_i)_{i=0}^n$ be an extended partition
on $(\Omega^{\calE,d},{\cal F}^{\sE,d})$.
Let $Y$ be a continuous $\R^{d'}$-valued process $(d'$ may
be different from $d)$,
and suppose that $Y$ is a semimartingale with
characteristic pair
$(B,C)$ under $\P$. If $\Delta(Y,T_i)$, $\Delta(B,T_i)$
and $\Delta(C,T_i)$ are all $\sG_i\vee\sigma(\Delta(X,T_i))$-measurable
for $i=1,\ldots,n$, then under $\P^{\otimes\Pi}$ the process
$Y$ is still a semimartingale with characteristic pair
$(B,C)$.
\end{proposition}

The proof of Proposition~\ref{P3.12} depends on
some preliminary results.\vadjust{\goodbreak}

%
%
\begin{lemma}\label{L3.14}
Let $\P_1$ and $\P_2$ be probability measures on
$(\Omega^{\calE,d},{\cal F}^{\sE,d})$, and let $T$ be a finite
(for every $\omega$)
$\{{\cal F}^{\sE,d}_t\}_{t\geq0}$-stopping time.
Let $M$ be a continuous
local martingale relative to $\{{\cal F}^{\sE,d}_t\}_{t\geq0}$
under $\P_1$ and $\P_2$.
Let $\sG$ be a sub-$\sigma$-field of ${\cal F}^{\sE,d}_T$ such that
$\P_1|_{\sG}\ll\P_2|_{\sG}$ and assume that
$\Mhat\define\Delta(M,T)$ is
$\sG\vee\sigma(\Delta(X,T))$-measurable.
Then $(M_t,{\cal F}^{\sE,d}_t)_{t\geq0}$
is a continuous local martingale under
$\P_{12}\define\P_1\otimes_{T,\sG}\P_2$.
\end{lemma}

\begin{pf}
It is sufficient to show that $M^T$ and $M - M^T$ are both
$\P_{12}$-local martingales relative to $\{{\cal F}^{\sE,d}_t\}_{t\geq
0}$. As
$M^T$ is a $\P_1$-local martingale, $M^T$ is ${\cal F}^{\sE
,d}_T$-measurable and
$\P_1$ and $\P_{12}$ agree on ${\cal F}^{\sE,d}_T$, we may
immediately conclude
that $M^T$ is a $\P_{12}$-local martingale.

For each
integer $n > 0$, define the stopping time
$S^n\define\inf\{t\geq T\dvtx|M_t-M_T|\geq n\}$. Then $M^{S_n} -
(M^{S_n})^T = M^{S_n} - M^T$ is bounded, and $\Delta(M^{S_n}, T)$ is
$\sG\vee\sigma(\Delta(X,T))$-measurable. As a result, we may assume
without loss of generality that $M - M^T$ is a uniformly integrable
$\P_{2}$-martingale.

We now show that $\Mhat$ is a $\P_{12}$-martingale with respect to the
filtration $\sFhat_t = \sF_{T+t}$. The process $\Mhat$ is clearly
$\{\sFhat_t\}_{t\geq0}$-adapted, and it follows from the
optional\vspace*{1pt}
sampling theorem that $\Mhat$ is a $\P_2$-martingale with respect to
the filtration $\{\sFhat_t\}_{t\geq0}$. For $0\leq s\leq t$, $A\in{\cal
F}^{\sE,d}_T$ and $B\in\sigma(\Delta_r(X,T)\dvtx0 \leq r \leq s)$, we
have from Corollary~\ref{C3.6}(b) that
\begin{eqnarray*}
\E_{12} \bigl[\ind_{A\cap
B} (\Mhat_t -
\Mhat_s ) \bigr] &=& \E_{12} \bigl[\ind_A
\E_{12} \bigl[\ind_B (\Mhat_t -
\Mhat_s ) | {\cal F}^{\sE,d}_T \bigr] \bigr]
\\
&=& \E_{12} \bigl[\ind_A \E_2 \bigl[
\ind_B (\Mhat_t -\Mhat_s ) | \sG\bigr]
\bigr]
\\
&=& \E_{12} \bigl[\ind_A\E_2 \bigl[
\ind_B \E_2 [\Mhat_t -\Mhat_s |
\sFhat_s ] | \sG\bigr] \bigr] = 0,
\end{eqnarray*}
where we have used the fact that $\Mhat$ is $\P_2$-martingale in
the last step. Writing $X^{T+s}_t = X^T_t + \Delta_{0\vee(t-T)
\wedge s}(X, T)$, we see that sets of the form $A\cap B$ generate
$\sFhat_s$. It then follows from Dynkin's $\pi$--$\lambda$ theorem
that $\Mhat$ is $\P_{12}$-martingale relative to $\{\sFhat_t\}_{t\geq
0}$.

To conclude the proof, we observe that
$(r - T)^+$ is a bounded $\{\sFhat_t\}_{t\geq0}$-stopping time and
$\sF_r \subset\sFhat_{(r-T)\vee0}$ for each $r \geq0$. Fixing $0
\leq s < t$ and $A \in\sF_s$, we have
\begin{eqnarray*}
\E_{12}\bigl[\ind_A \bigl(M_t -
M^T_t\bigr) \bigr] &=& \E_{12}[
\ind_A \Mhat_{(t-T)\vee0}] = \E_{12}[
\ind_A \Mhat_{(s-T)\vee0}] \\
&=& \E_{12}\bigl[
\ind_A \bigl(M_s - M^T_s\bigr)
\bigr],
\end{eqnarray*}
so $M - M^T$ is $\P_{12}$-martingale relative to
$\{\sF_t\}_{t\geq0}$.
\end{pf}

%
%
\begin{lemma}\label{L3.15}
Let $\P$ be a probability measure on $(\Omega^{\calE,d},{\cal
F}^{\sE,d})$,
and let $M$ be a uniformly integrable $\P$-martingale
relative to $\{{\cal F}^{\sE,d}_t\}_{t\geq0}$. Let $S$, $T$ and $U$
be stopping times with $T\leq U$ almost surely, and let
$Z$ be an ${\cal F}^{\sE,d}_T$-measurable bounded random variable. Then
$\E[(M_U-M_T)Z | {\cal F}^{\sE,d}_S ]
=(M_{U\wedge S}-M_{T\wedge S})Z$.
\end{lemma}

\begin{pf}
Because
\[
\ind_{\{S\leq T\}}\E\bigl[(M_U-M_T)Z | {\cal
F}^{\sE,d}_S \bigr] =\ind_{\{S\leq T\}}\E\bigl[Z\E
\bigl[(M_U-M_T) | {\cal F}^{\sE,d}_T
\bigr] | {\cal F}^{\sE,d}_S \bigr] =0,
\]
we have
\begin{eqnarray*}
&&\E\bigl[(M_U-M_T)Z | {\cal
F}^{\sE,d}_S \bigr]
\\[-2pt]
&&\qquad= \ind_{\{T<S\leq U\}}\E\bigl[(M_U-M_T)Z | {\cal
F}^{\sE,d}_S \bigr] +\ind_{\{U<S\}}\E
\bigl[(M_U-M_T)Z | {\cal F}^{\sE,d}_S
\bigr]
\\[-2pt]
&&\qquad= \ind_{\{T<S\leq U\}}(M_S-M_T)Z +
\ind_{\{U<S\}}(M_U-M_T)Z
\\[-2pt]
&&\qquad= (M_{U\wedge S}-M_{T\wedge S})Z.
\end{eqnarray*}
\upqed
\end{pf}

%
%
\begin{lemma}\label{L3.16}
Let $\P_1$ and $\P_2$
be probability measures on $(\Omega^{\calE,d},{\cal F}^{\sE,d})$
and let $T$ be a finite (for every $\omega$)
$\{{\cal F}^{\sE,d}_t\}_{t\geq0}$-stopping time.
Let $M^1$, $M^2$ and $C$ be continuous
$\{{\cal F}^{\sE,d}_t\}_{t\geq0}$-adapted real-valued processes
such that $M^1$, $M^2$ and $M^3\define M^1M^2-C$ are local
martingales relative to $\{{\cal F}^{\sE,d}_t\}_{t\geq0}$ under $\P_1$
and $\P_2$.
Let $\sG$ be a sub-$\sigma$-field of ${\cal F}^{\sE,d}_T$ such that
$\P_1|_{\sG}\ll\P_2|_{\sG}$ and assume that
$\Mhat^1\define\Delta(M^1,T)$,
$\Mhat^2\define\Delta(M^2,T)$ and $\Chat\define\Delta(C,T)$
are $\sG\vee\sigma(\Delta(X,T))$-measurable.
Then $(M^3,{\cal F}^{\sE,d}_t)_{t\geq0}$ is a
local martingale under \mbox{$\P_{12}\define\P_1\otimes_{\sG,T}\P_2$}.
\end{lemma}

\begin{pf}
We cannot apply Lemma~\ref{L3.14} directly because we
did not assume that
$\Mhat_3$ is $\sG\vee\sigma(\Delta(X,T))$-measurable.
Instead, we define the process
\begin{eqnarray*}
Y_t &\define& M_{t\wedge T}^1M_{t\wedge T}^2
+\bigl(M_t^1-M_{t\wedge T}^1\bigr)
\bigl(M_t^2-M_{t\wedge T}^2
\bigr)-C_t
\\[-2pt]
&=& M_t^3-\bigl(M_t^1-M_{t\wedge T}^1
\bigr)M_{t\wedge T}^2 -\bigl(M_t^2-M_{t\wedge T}^2
\bigr)M_{t\wedge T}^1,\qquad t\geq0,
\end{eqnarray*}
for which $\Delta(Y,T)=\Delta(M^1,T)\Delta(M^2,T)-\Delta(C,T)$ is
$\sG\vee\sigma(\Delta(X,T))$-measur\-able.
Define
$T_n\define
\inf\{t\geq0\dvtx|M_t^1|\vee|M_t^2|\vee|M_t^3|
\vee|C_t|\geq n \}$,
and set $M^{i,n}\define(M^i)^{T_n}$ for $i=1,2,3$,
$C^n\define C^{T_n}$ and $Y^n\define Y^{T_n}$. For fixed $n$,
the processes $M^{i,n}$, $i=1,2,3$, and $Y^n$ are bounded.
For $0\leq s\leq t$,
we apply Lemma~\ref{L3.15} with $M=M^{1,n}$, $Z=M_{t\wedge T}^{2,n}$,
$S=s$, $T=t\wedge T$ and $U=t$, and use the fact that
$M_s^{1,n}-M_{s\wedge T}^{1,n}=0$ if $T\geq s$ to obtain
%
%
\begin{eqnarray} \label{3.24}\qquad
\E_k \bigl[\bigl(M_t^{1,n}-M_{t\wedge T}^{1,n}
\bigr)M_{t\wedge T}^{2,n} | {\cal F}^{\sE,d}_s
\bigr] &=& \bigl(M_s^{1,n}-M_{s\wedge T}^{1,n}
\bigr)M_{t\wedge T}^{2,n}
\nonumber\\[-9pt]\\[-9pt]
&=& \bigl(M_s^{1,n}-M_{s\wedge T}^{1,n}
\bigr)M_{s\wedge T}^{2,n},\qquad k=1,2.\nonumber
\end{eqnarray}
The same equality holds if we reverse the roles of $M^{1,n}$
and $M^{2,n}$. Finally, because $M^{3,n}$ is a martingale,
\begin{eqnarray*}
\E_k\bigl[Y_t^n | {\cal
F}^{\sE,d}_s\bigr] &=& \E_k\bigl[M_t^{3,n}
| {\cal F}^{\sE,d}_s\bigr] -\E_k \bigl[
\bigl(M_t^{1,n}-M_{t\wedge T}^{1,n}
\bigr)M_{t\wedge T}^{2,n} | {\cal F}^{\sE,d}_s
\bigr]
\\[-2pt]
&&{} -\E_k \bigl[\bigl(M_t^{2,n}-M_{t\wedge T}^{2,n}
\bigr)M_{t\wedge T}^{1,n} | {\cal F}^{\sE,d}_s
\bigr]
\\
&=& M_s^{3,n}-\bigl(M_s^{1,n}-M_{s\wedge T}^{1,n}
\bigr)M_{s\wedge T}^{2,n} -\bigl(M_s^{2,n}-M_{s\wedge T}^{2,n}
\bigr)M_{s\wedge T}^{1,n}
\\
&=& Y_s^n,\qquad k=1,2,
\end{eqnarray*}
so $Y$ is a local martingale under both $\P_1$ and $\P_2$.
Lemma~\ref{L3.14} implies that $M^1$, $M^2$ and $Y$
are $\P_{12}$-local martingales. Therefore, (\ref{3.24})
holds under $\P_{12}$ as well,
from which we conclude that
\begin{eqnarray*}
\E_{12}\bigl[M_t^{3,n} | {\cal F}^{\sE,d}_s\bigr]
&=& \E_{12}
\bigl[Y_t^n | {\cal F}^{\sE,d}_s
\bigr] +\E_{12} \bigl[\bigl(M_t^{1,n}-M_{t\wedge T}^{1,n}
\bigr)M_{t\wedge T}^{2,n} | {\cal F}^{\sE,d}_s
\bigr]
\\
&&{} +\E_{12} \bigl[\bigl(M_t^{2,n}-M_{t\wedge T}^{2,n}
\bigr)M_{t\wedge T}^{1,n} | {\cal F}^{\sE,d}_s
\bigr]
\\
&=& Y_s^n+\bigl(M_s^{1,n}-M_{s\wedge T}^{1,n}
\bigr)M_{s\wedge T}^{2,n} +\bigl(M_s^{2,n}-M_{s\wedge T}^{2,n}
\bigr)M_{s\wedge T}^{1,n}
\\
&=& M_s^{3,n},\qquad 0\leq s\leq t.
\end{eqnarray*}
\upqed
\end{pf}

\begin{pf*}{Proof of Proposition~\ref{P3.12}}
According to Remark~\ref{R3.8}, $\P^{\otimes\Pi}=\P^{n+1}$,
where $\P^i$ is defined
recursively by $\P^0=\P$ and
$\P^{i+1}=\P^i\otimes_{T_i,\sG_i}\P$, $i=0,\ldots,n$.
If $M$ is a continuous local martingale under $\P$
and $\Delta(M,T_i)$ is
$\sG_i\vee\sigma(\Delta(X,T_i))$-measurable
for $i=0,\ldots,n$, then
repeated application of Lemma~\ref{L3.14} shows that
$M$ is a $\P^i$-local martingale for
$i=1,\ldots,n,n+1$, and in particular,
$M$ is a continuous local martingale
under $\P^{\otimes\Pi}$. Similarly, if $M^1$ and $M^2$
are continuous local martingales under $\P$,
$C$ is a finite variation process such that
$M^3\define M^1M^2-C$ is a local martingale under $\P$,
and $\Delta(M^1,T_i)$, $\Delta(M^2,T_i)$ and $\Delta(C,T_i)$
are $\sG_i\vee\sigma(\Delta(X,T_i))$-measurable for $0=1,\ldots,n$,
then repeated application of Lemma~\ref{L3.16} shows
that $M^1M^2-C$ is a $\P^i$-local
martingale for $i=1,\ldots,n,n+1$.
In particular, $M^1M^2-C$ is a continuous local martingale
under $\P^{\otimes\Pi}$.
These observations combined
with Proposition~\ref{P3.9}(a)
prove the desired result.
\end{pf*}

\section{Conditional expectations}
\label{SecCondExp}

The results of this section are implicit in
Krylov~\cite{Kry3} and Gy\"ongy~\cite{G}.
We use the notation introduced in Sections
\ref{SecMainResult} and~\ref{SConMea}. In addition,
we denote the Borel $\sigma$-field on $[0,t]$
by $\sB[0,t]$ and the Borel $\sigma$-field on
$[0,\infty)$ by $\sB[0,\infty)$.

%
%
\begin{proposition}\label{P4.1}
Let $Z$ be an $\sE$-valued process and let
$\Gamma$ be an $\R^d$-valued process
(resp., a $d\times d$ matrix-valued process) taking
values in a closed convex set $K$, and satisfying
$\E[\int_0^t\|\Gamma_u\| \,du]
<\infty$ for all $t\geq0$.
Then there
exists an $\R^d$-valued measurable
function (resp., a $d\times d$ matrix-valued
measurable function) $\Gammahat$, defined on
$[0,\infty)\times\sE$, taking values in $K$,
and there exists a Lebesgue-null set $N\subset[0,\infty)$,
so that
%
%
\begin{equation}
\label{4.2} \Gammahat(t,Z_t)=\E[\Gamma_t |
Z_t],\qquad \P\mbox{-a.s.}, t\in N^c.
\end{equation}
\end{proposition}

The proof of Proposition~\ref{P4.1} depends on the following
lemma.

%
%
\begin{lemma}\label{L4.2}
Let $Z$ be an $\sE$-valued process and let
$\Gamma$ be a real-valued process satisfying
$\E\int_0^t|\Gamma_u| \,du<\infty$ for all $t\geq0$.
Let $\Gammahat$ be a real-valued measurable function
on $[0,\infty)\times\sE$. There exists a
Lebesgue-null set $N\subset[0,\infty)$ so that (\ref{4.2})
holds if and only if for every bounded
$\sB[0,\infty)\otimes\frakE$-measurable real-valued
function $f$,
%
%
\begin{equation}
\label{4.4} \E\int_0^t\Gammahat(u,Z_u)f(u,Z_u)
\,du =\E\int_0^t\Gamma_uf(u,Z_u)
\,du,\qquad t\geq0.
\end{equation}
\end{lemma}

\begin{pf} If (\ref{4.2}) holds, then (\ref{4.4}) follows
from Fubini's theorem.

To prove the converse, we assume (\ref{4.4}).
Taking $f(u,Z_u)=\mbox{sgn}(\Gammahat(u,Z_u))$
and using the integrability of $\Gamma$, we
see that $\E\int_0^t|\Gammahat(u,Z_u)| \,du<\infty$
for all $t\geq0$.

The $\sigma$-field $\frakE$ is generated by
a collection of open balls intersected
with~$\sE$, each ball having a rational radius and centered
at a point in a countable dense subset of the
separable metric space containing $\sE$. Let $\sO$
denote the collection of finite intersections of
this countable collection of sets. Then $\sO$
is itself countable and $\frakE=\sigma(\sO)$.
We enumerate the sets in $\sO$ as
$O_1,O_2,\ldots\,$.
Define
$g_n(t)\define
\E[ (\Gammahat(t,Z_t)-\Gamma_t )
\ind_{\{Z_t\in O_n\}} ]$.
For $B\in\sB[0,t]$, (\ref{4.4}) implies
\[
\int_Bg_n(u) \,du= \E\int_0^t
\bigl(\Gammahat(u,Z_u)-\Gamma_u \bigr)
\ind_{\{(u,Z_u)\in B\times O_n\}} \,du=0.
\]
Since both $t\geq0$ and $B\in\sB[0,t]$ are arbitrary,
we conclude that $g_n=0$ for Lebesgue-almost every $t\geq0$.
Thus,
$N\define\{t\geq0|g_n(t)\neq0\mbox{ for some }n\}$
is a Lebesgue-null set.

The collection of sets $A\in\frakE$ for which
%
%
\begin{equation}
\label{4.5} \E\bigl[\bigl(\Gammahat(t,Z_t)-\Gamma_t\bigr)
\ind_{\{Z_t\in A\}}\bigr]=0,\qquad t\in N^c,
\end{equation}
is a $\lambda$ system containing $\sO$,
and the Dynkin $\pi$--$\lambda$
theorem implies that (\ref{4.5}) holds for
every $A\in\frakE$. This gives us (\ref{4.2}).
\end{pf}

\begin{pf*}{Proof of Proposition~\ref{P4.1}}
Except for the assertion that $\Gammahat$ takes values
in the set $K$, it suffices to prove
the proposition for the case that $\Gamma$
is real-valued. We can then apply the one-dimensional
result to each
component of the $\Gamma$ in the proposition.

In the one-dimensional case, we define the $\sigma$-finite
measure
\[
\mu(A)\define\E\int_0^{\infty}\ind_A(u,Z_u)
\,du,\qquad A\in\sB[0,\infty)\otimes\frakE,
\]
and the $\sigma$-finite signed measure
\[
\nu(A)\define\E\int_0^{\infty}\Gamma_u
\ind_A(u,Z_u) \,du,\qquad A\in\sB[0,\infty)\otimes\frakE.
\]

Obviously, $\nu\ll\mu$, so we can define
$\Gammahat(t,z)=\frac{d\nu}{d\mu}(t,z)$ for
$(t,z)\in[0,\infty]\otimes\sE$.
Let $f$ be a bounded $\sB[0,\infty)\otimes\frakE$-measurable
real-valued function. For $t\geq0$,
\begin{eqnarray*}
\E\int_0^t\Gammahat(u,Z_u)f(u,Z_u)
\,du &=& \int_{[0,t]\times\sE} \Gammahat(u,z)f(u,z)\mu(du,dz)
\\
&=& \int_{[0,t]\times\sE} f(u,z) \nu(du,dz)
\\
&=& \E\int_0^t\Gamma_uf(u,Z_u)
\,du.
\end{eqnarray*}
Equation (\ref{4.2}) follows from Lemma~\ref{L4.2}.

Let us now consider the case of a multi-dimensional $\Gamma$
taking values in a closed convex set $K$.
We have already shown the existence of $\Gammahat$ such
that (\ref{4.2}) holds, and it remains to show
that $\Gammahat$ takes values in $K$.
Define $\varphi\dvtx\R^d\rightarrow\R$
(resp., $\varphi\dvtx\R^d\times\R^d\rightarrow\R$) by
$\varphi(\gamma)={\min_{\kappa\in K}}\|\gamma-\kappa\|$,
which is the distance from $\gamma$ to $K$. One can verify
from the triangle inequality that for each constant $c$,
the set $\{\gamma|\varphi(\gamma)\leq c\}$ is convex,
and hence, $\varphi$ is a continuous convex real-valued function.
Such a function has the property that
$\varphi(\gamma)=\max\{\ell(\gamma)|
\ell\mbox{ is linear and }\ell\leq\varphi\}$.
This permits us to establish the Jensen inequality
\begin{eqnarray*}
\E\bigl[\varphi(\Gamma_t) | Z_t \bigr] &\geq& \max
\bigl\{\E\bigl[\ell(\Gamma_t) | Z_t \bigr] | \ell\mbox{
is linear and }\ell\leq\varphi\bigr\}
\\
&=& \max\bigl\{\ell\bigl(\E[\Gamma_t | Z_t] \bigr) |
\ell\mbox{ is linear and }\ell\leq\varphi\bigr\}
\\
&=& \varphi\bigl(\E[\Gamma_t | Z_t] \bigr)
\\
&=& \varphi\bigl(\Gammahat(t,Z_t) \bigr),\qquad t\in N^c.
\end{eqnarray*}
But $\Gamma$ takes values in $K$,
so the left-hand side of this inequality is zero.
Thus the right-hand side is zero,
implying $\Gammahat(t,Z_t)\in K$ almost surely
for each $t\in N^c$.
We can modify $\Gammahat(t,z)$
so that it takes values in $K$ for every
$t$, and (\ref{4.2}) still holds.
\end{pf*}

%
%
\begin{definition}\label{D4.3}
Let $\{\Gamma^i\}_i$ be a collection of processes on
some probability space $(\Omega,\sF,\P)$
and let $T$ be a $[0,\infty)$-valued random
variable. We say the collection $\{\Gamma^i\}_i$
is \textit{strongly independent} of
$T$ if there is a $\sigma$-field $\sG\subset\sF$
such that each $\Gamma_i$ is $\sB[0,\infty)\times\sG$-measurable
and $\sG$ is independent of $\sigma(T)$.
\end{definition}

%
%
\begin{proposition}\label{P4.5}
Within the setting of Proposition~\ref{P4.1},
let $T$ be a $[0,\infty)$-valued random variable
whose distribution $\mu\define\P\circ T^{-1}$
is absolutely continuous with respect to Lebesgue
measure. Assume also that the pair
of processes $(\Gamma, Z)$ is strongly
independent of $T$ and $\E|\Gamma_T|<\infty$. Then
%
%
\begin{equation}
\label{4.7} \Gammahat(T,Z_T)=\E[\Gamma_T |
T,Z_T],\qquad \P\mbox{-a.s.}
\end{equation}
\end{proposition}

\begin{pf}
We first observe that
%
%
\begin{equation}
\label{4.6} \int_0^{\infty}\E[\Xi_t]
\mu(dt) =\E[\Xi_T]\vadjust{\goodbreak}
\end{equation}
for any process $\Xi$ that is strongly independent
of $T$ and
satisfies $\E|\Xi_T|<\infty$.
To see this, consider the case
$\Xi_t=\sum_{i=1}^n\ind_{A_i}\ind_{B_i}(t)$,
where $A_i\in\sG$, the $\sigma$-field
in Definition~\ref{D4.3} and $B_i\in\sB[0,\infty)$.
Then use the monotone class theorem.

Now let $f\dvtx[0,\infty)\times\sE\rightarrow\R$
be a bounded, $\sB[0,\infty)\otimes\frakE$-measurable
real-valued function. Proposition~\ref{P4.1}
implies
$\E[\Gammahat(t,Z_t)f(t,Z_t) ]
=\E[\Gamma_tf(t,Z_t) ]$ for all \mbox{$t\in N^c$}.
Integrating both sides of this equation
with respect to $\mu(dt)$
and using (\ref{4.6}), we obtain
$\E[\Gammahat(T,Z_T)f(T,Z_T) ]
=\E[\Gamma_Tf(T,Z_T) ]$.
Equation (\ref{4.7}) follows.
\end{pf}

\section{Approximation}\label{S5}

We collect in this section three approximation results needed to prove
Theorem~\ref{T2.4}. We denote by $\N$ the set of natural numbers and
define $\Nbar\define\N\cup\{\infty\}$. We recall that
$\lambda_{[0,t]}$
denotes Lebesgue measure on $[0,t]$.

\subsection{Convergence of the integral
of a process}\label{S5.1}

%
%
\begin{proposition}\label{P5.2}
Let $\{Z^m\}_{m\in\Nbar}$ be a collection
of continuous $\sE$-valued processes,
possibly defined on different probability
spaces under different probability
measures $\Q^m$.
Let
$f\dvtx[0,\infty)\times\sE\rightarrow\R^d$
be a measurable function. Assume:
\begin{longlist}[(vii)]
\item[(i)]
for each $t\in[0,\infty)$,
the distribution of $Z_t^m$ under $\Q^m$ is independent of
$m\in\N$,
\item[(ii)]
the distribution on $C^{\sE}$
of $Z^m$ under $\Q^m$ converges weakly to the distribution
of $Z^{\infty}$ under $\Q^{\infty}$, that is,
$\Q^m\circ(Z^m)^{-1}\Rightarrow\Q^{\infty}\circ(Z^{\infty})^{-1}$
and
\item[(iii)]
$\E^{\Q^1}\int_0^t\|f(u,Z_u^1)\| \,du
<\infty$ for every $t\in[0,\infty)$.
\end{longlist}
Then:
\begin{longlist}[(vii)]
\item[(iv)]
for every $m\in\Nbar$ the
integral process
$F_t^m\define\int_0^tf(s,Z_s^m) \,ds$,
$t\in[0,\infty)$,
is defined $\Q^m$-almost surely,
\item[(v)]
$\Q^m[F^m\in C^d]=1$ for every
$m\in\Nbar$,
\item[(vi)]
$\{f(\cdot,Z^m_{\cdot}),\lambda_{[0,t]}\times\Q^m
\}_{m\in\Nbar}$ is uniformly integrable
for every $t\in[0,\infty)$,
\item[(vii)]
$(Z^m,F^m)\Rightarrow(Z^{\infty},F^{\infty})$.
\end{longlist}
\end{proposition}

\begin{pf}
It suffices to prove parts (iv)--(vi)
of the lemma for
the case $d=1$, since these results
can be applied component-wise
to the $d$-dimensional~$f$.

Define the measure $\mu$ on $[0,\infty)\times\sE$ by
$\mu(A)\define\E^{\Q^m} \int_0^{\infty}\ind_A(s,Z^m_s) \,ds$.
Assumption (i) and the convergence in (ii) imply that the distribution
of $Z^m_t$ is independent of $m\in\Nbar$,\vadjust{\goodbreak} so it does not matter which
$m\in\Nbar$ we use in the definition of $\mu$.
Therefore, for each $m\in\Nbar$ and $M>0$,
\begin{eqnarray*}
\E^{\Q^m}\int_0^t\bigl|f
\bigl(s,Z^m_s\bigr)\bigr| \ind_{\{|f(x,Z^m_s)|\geq M\}} \,ds &=& \int
_{[0,t]\times\sE} \bigl|f(s,e)\bigr|\ind_{\{|f(s,e)|\geq M\}} \mu(ds,de)
\\
&=& \E^{\Q^1}\int_0^t\bigl|f
\bigl(s,Z_s^1\bigr)\bigr| \ind_{\{|f(s,Z_s^1)|\geq M\}} \,ds.
\end{eqnarray*}
Setting $M=0$, we obtain (iv) and (v) from (iii).
Condition (iii) implies that the
last term can be made arbitrarily
small by choosing $M$ large,
and (vi) also follows.

To prove (v), we use (vi) and
Lusin's theorem
to choose for each $k\in\N$ a bounded
continuous function
$f^k\dvtx[0,k]\times\sE\rightarrow\R^d$ such
that
\[
\lim_{k\rightarrow\infty} \int_{[0,k]\times\sE} \bigl\|f(t,e)-f^k(t,e)
\bigr\| \mu(dt,de)=0.
\]
The mapping
$z\mapsto\int_0^{\cdot\wedge k}f^k(s,z(s))\,ds$
is continuous from $C^{\sE}$ to
$C^d$, which implies that
%
%
\begin{equation}
\label{5.8} \bigl(Z^m,F^{m,k}\bigr)\Rightarrow
\bigl(Z^{\infty},F^{\infty,k}\bigr)\qquad\mbox{as }m\rightarrow\infty,
\end{equation}
where
$F^{m,k}_t\define
\int_0^{t\wedge k}f^k(s,Z_s^m)\,ds$.
But for each fixed $T$ and $k\geq T$,
%
%
\begin{eqnarray}\label{5.8a}
&&
\sup_{m\in\Nbar} \Q^m \Bigl[\sup_{0\leq t\leq T}
\bigl\|F_t^m-F_t^{m,k}\bigr\|>\varepsilon\Bigr]\nonumber\\
&&\qquad\leq \sup_{m\in\Nbar} \frac{1}{\varepsilon} \E^{\Q^m}
\sup_{0\leq t\leq T} \bigl\|F^m_t-F^{m,k}_t
\bigr\|
\\
&&\qquad\leq \frac{1}{\varepsilon} \int_{[0,k]\times\sE} \bigl\|f(s,e)-f^k(s,e)
\bigr\| \mu(ds,de),\nonumber
\end{eqnarray}
which has limit zero as $k\rightarrow\infty$.
In particular, the convergence
$F^{m,k}\Rightarrow F^m$ as
$k\rightarrow\infty$ is uniform
in $m\in\Nbar$.

Let $\Psi\dvtx C^{\sE}\times C^1\rightarrow
\R$ be a uniformly continuous bounded function.
To prove weak convergence of measures on a metric space,
it suffices to consider such functions (see~\cite{P},
Chapter II, Theorem 6.1).
We have
\begin{eqnarray*}
&&
\bigl|\E^m \bigl[\Psi\bigl(Z^m,F^m
\bigr)\bigr] -\E^\infty\bigl[\Psi\bigl(Z^{\infty},F^{\infty}
\bigr) \bigr] \bigr|
\\
&&\qquad= \bigl|\E^m \bigl[\Psi\bigl(Z^m,F^m\bigr)
-\Psi\bigl(Z^m,F^{m,k}\bigr) \bigr] \bigr| \\
&&\qquad\quad{}+ \bigl|\E^m
\bigl[\Psi\bigl(Z^m,F^{m,k}\bigr)\bigr] -\E^\infty
\bigl[\Psi\bigl(Z^{\infty},F^{\infty,k}\bigr) \bigr] \bigr|
\\
&&\qquad\quad{} + \bigl|\E^\infty\bigl[\Psi\bigl(Z^{\infty},F^{\infty,k}\bigr)
-\Psi\bigl(Z^{\infty},F^{\infty}\bigr) \bigr] \bigr|.
\end{eqnarray*}
Given $\varepsilon>0$, (\ref{5.8a}) guarantees that we can choose $k$
so large that the first and third terms on the right-hand side are less
than $\varepsilon$, independently of $m$. For this value of~$k$, we can
then use (\ref{5.8}) to choose $M$ so that for all $m\geq M$, the
second term is also less than $\varepsilon$.
\end{pf}

\subsection{Approximation by step functions}
\label{S5.2}

We show in Proposition~\ref{P5.5} below
that an arbitrary integrable process
can be approximated in $L_1(\P\times\lambda_{[0,t]})$
by step functions obtained by sampling
the process at random partition points.

%
%
\begin{lemma}\label{L5.4}
Let $f\dvtx[0,\infty)\rightarrow\R^d$
be a measurable function with\break
$\int_0^t\|f(s)\|\,ds<\infty$
for every $t\in[0,\infty)$.
Define the sets
%
%
\begin{equation}
\label{5.9a} I_i^n\define\biggl\{(t,u)\in[0,\infty)
\times[0,1]\dvtx\frac{u+i-1}{n}\leq t<\frac{u+i}{n} \biggr\},\qquad
i=1,2,\ldots,\hspace*{-28pt}
\end{equation}
and define the sequence of
functions
$f_n(t,u)=\sum_{i=1}^{\infty}
f (\frac{u+i-1}{n} )
\ind_{I_i^n}(t,u)$.
Then
\[
\lim_{n\rightarrow\infty} \int_0^1\int
_0^t \bigl\|f(s)-f_n(s,u)\bigr\| \,ds \,du=0,\qquad
t\in[0,\infty).
\]
\end{lemma}

\begin{pf}
Fix $t>0$ and $\varepsilon>0$.
Choose a continuous, $\R^d$-valued
function $g$ defined on $[0,t+1]$
for which
$\int_0^{t+1}\|f(s)-g(s)\|\,ds\leq\varepsilon$.
Set $m=\lceil t \rceil\in[t,t+1)\cap\N$
and set
$
g_n(s,u)\define\sum_{i=1}^{mn}
g (\frac{u+i-1}{n} )
\ind_{I_i^n}(s,u).
$
We have
\begin{eqnarray*}
&&\int_0^1\int_0^t
\bigl\|f_n(s,u)-g_n(s,u)\bigr\| \,ds \,du
\\
&&\qquad\leq\sum_{i=1}^{mn}\int
_0^1 \int_{({u+i-1})/{n}}^{({u+i})/{n}}
\biggl\llVert f \biggl(\frac{u+i-1}{n} \biggr) -g \biggl(\frac{u+i-1}{n}
\biggr)\biggr\rrVert\,ds \,du
\\
&&\qquad= \sum_{i=1}^{mn}\int
_0^1\biggl\llVert f \biggl(\frac{u+i-1}{n}
\biggr) -g \biggl(\frac{u+i-1}{n} \biggr)\biggr\rrVert\,\frac{du}{n}
\\
&&\qquad= \sum_{i=1}^{mn}\int
_{({i-1})/{n}}^{{i}/{n}} \bigl\|f(v)-g(v)\bigr\| \,dv \leq\varepsilon.
\end{eqnarray*}
Because $g$ is uniformly continuous on $[0,t+1]$,
we may choose $N$ so that $\|g(s_2)-g(s_1)\|\leq\varepsilon/t$
whenever $|s_2-s_1|\leq1/N$. By enlarging $N$
if necessary, we can also ensure that
$\int_0^{1/N}\|g(s)\|\,ds\leq\varepsilon$.
Therefore, for $n\geq N$, we have
\begin{eqnarray*}
&&
\int_0^1\int
_0^t\bigl\|f(s)-f_n(s,u)\bigr\| \,ds \,du
\\
&&\qquad\leq\int_0^t\bigl\|f(s)-g(s)\bigr\| \,ds +\int
_0^1\int_{u/n}^t
\bigl\|g(s)-g_n(s,u)\bigr\| \,ds \,du
\\
&&\qquad\quad{} +\int_0^{1/n}\bigl\|g(s)\bigr\| \,ds +\int
_0^1\int_0^t
\bigl\|g_n(s,u)-f_n(s,u)\bigr\| \,ds \,du \leq4\varepsilon.
\end{eqnarray*}
\upqed
\end{pf}

%
%
\begin{proposition}\label{P5.5}
Let $(\Omega',\sF',\Q')$ be a probability space that
supports an $\R^d$-valued process $a$ satisfying
%
%
\begin{equation}
\label{5.10} \E^{\Q'}\int_0^t
\|a_s\| \,ds<\infty,\qquad t\geq0.
\end{equation}
Set $\Omega\define[0,1]\times\Omega'$, with generic
point $\omega=(u,\omega')$, and define
$U(u,\omega')=u$. Set $\sF=\sB[0,1]\otimes\sF'$,
$\Q=\lambda_{[0,1]}\times\Q'$ and extend $a$ to
$\Omega$ via the abuse of notation
$a(u,\omega')\define a(\omega')$.
Finally, define the random times $T_0^n\define0$,
$T_i^n\define(U+i-1)/n$ for $i=1,2,\ldots,n^2$
and $T_{n^2+1}^n\define\infty$. Then the
\textit{sampled process}
\[
a^n_t(\omega) \define\sum
_{i=1}^{n^2}a_{T_i^n(\omega)}(\omega)
\ind_{[T_i^n(\omega),T_i^{n+1}(\omega))}(t)= \sum_{i=1}^{n^2}a_{(u+i-1)/n}
\bigl(\omega'\bigr) \ind_{I_i^n}(t,u),
\]
where
$I_i^n$ is defined by (\ref{5.9a}),
satisfies
%
%
\begin{equation}
\label{5.11} \lim_{n\rightarrow\infty}\E^{\Q}\int_0^t
\bigl\|a_s-a_s^n\bigr\| \,ds=0,\qquad t\geq0.
\end{equation}
\end{proposition}

\begin{pf}
Define
$A^n_t(\omega')\define\int_0^1\int_0^t
\|a_s(\omega')-a_s^n(u,\omega')\| \,ds \,du$ for $t\geq0$
and \mbox{$\omega'\in\Omega'$}.
Assumption (\ref{5.10}) implies that
$\int_0^t\|a_s(\omega')\| \,ds<\infty$ for all
$t\geq0$ for $\Q'$-almost every
$\omega'$. For fixed $\omega'$ satisfying this condition,
Lemma~\ref{L5.4} then shows that
$\lim_{n\rightarrow\infty}A_t^n(\omega')=0$ for every $t\geq0$.
Equation (\ref{5.11}) is equivalent to
\[
\lim_{n\rightarrow0}\E^{\Q'}A_t^n=0,\qquad t\geq0,
\]
and to obtain this result it now suffices to show that
for each fixed $t\geq0$, the collection of random variables
$\{A_t^n\}_{n=1}^{\infty}$ is uniformly integrable under $\Q'$.

We first show that $\{a^n\}_{n=1}^{\infty}$ is uniformly
integrable with respect
to $\lambda_{[0,t]}\times\Q$ for every $t\geq0$.
Toward this end, fix $t\geq0$ and set $m\define\lceil t\rceil
\in[t,t+1)\cap\N$, so that $t\leq T_{mn+1}^n$.
Then
\begin{eqnarray*}
&&
\E^{\Q} \bigl[\|a_{T_i^n}\|\ind_{\{\|a_{T_i^n}\|\geq M\}} \bigr] \\
&&\qquad= \int
_0^1\E^{\Q'} \bigl[
\bigl\|a_{(u+i-1)/n}^n\bigr\| \ind_{\{\|a_{(u+i-1)/n}\|\geq M\}} \bigr]
\,\frac{du}{n}
\\
&&\qquad= \int_{(i-1)/n}^{i/n} \E^{\Q'} \bigl[
\|a_s\|\ind_{\{\|a_s\|\geq M\}} \bigr] \,ds,\qquad i=1, \ldots, mn,
\end{eqnarray*}
and
\begin{eqnarray*}
\E^{\Q}\int_0^t\bigl\|a_s^n
\bigr\|\ind_{\{\|a_s^n\|\geq M\}} \,ds &\leq& \E^{\Q}\int_0^{T_{mn+1}^n}
\bigl\|a_s^n\bigr\|\ind_{\{\|a_s^n\|\geq M\}} \,ds
\\
&=& \frac{1}{n}\sum_{i=1}^{mn}
\E^{\Q} \bigl[\|a_{T_i^n}\|\ind_{\{\|a_{T_i^n}\|\geq M\}} \bigr]
\\
&=& \E^{\Q'}\int_0^m
\|a_s\|\ind_{\{\|a_s\|\geq M\}} \,ds.
\end{eqnarray*}
The uniform integrability of $\{a^n\}_{n=1}^{\infty}$ under
$\lambda_{[0,t]}\times\Q$ follows from (\ref{5.10}).
This implies the uniform
integrability of $\{\|a-a^n\|\}_{n=1}^{\infty}$.
Jensen's inequality implies
\begin{eqnarray*}
\E^{\Q'} \bigl[\bigl(A_t^n-M\bigr)^+ \bigr]
&=& \E^{\Q'} \biggl[ \biggl(\int_0^1
\int_0^t \bigl\|a_s(
\cdot)-a_s^n(u,\cdot)\bigr\| \,ds \,du-M \biggr)^+ \biggr]
\\
&\leq& \E^{\Q'} \biggl[\int_0^1
\biggl(\int_0^t \bigl\|a_s(
\cdot)-a_s^n(u,\cdot)\bigr\| \,ds-M \biggr)^+\,du \biggr]
\\
&=& \E^{\Q} \biggl[ \biggl(\int_0^t
\bigl\|a_s-a_s^n\bigr\| \,ds-M \biggr)^+ \biggr]
\\
&\leq& \E^{\Q} \biggl[\int_0^t
\biggl(\bigl\|a_s-a_s^n\bigr\| -\frac{M}{t}
\biggr)^+ \,ds \biggr],
\end{eqnarray*}
and the uniform integrability of $\{\|a-a^n\|\}_{n=1}^{\infty}$ under
$\lambda_{[0,t]}\times\Q$ implies
that for every $\varepsilon>0$, there exists $M_{\varepsilon}>0$ such that
\[
\sup_{n\in\N}\E^{\Q'} \bigl[\bigl(A_t^n-M_{\varepsilon}
\bigr)^+ \bigr]\leq\varepsilon.
\]
Consequently,
\begin{eqnarray*}
\sup_{n\in\N}\E^{\Q'} \bigl[A_t^n
\ind_{\{A_t^n\geq2M_{\varepsilon}\}} \bigr] &=& \sup_{n\in\N} \E^{\Q'}
\bigl[
\bigl(A_t^n-2M_{\varepsilon}\bigr)^+
+2M_{\varepsilon}\ind_{\{A_t^n\geq2M_{\varepsilon}\}} \bigr]
\\
&\leq& \sup_{n\in\N} \E^{\Q'} \bigl[\bigl(A_t^n-2M_{\varepsilon}
\bigr)^+ +2\bigl(A_t^n-M_{\varepsilon}\bigr)^+ \bigr]
\\
&\leq& 3\sup_{n\in\N} \E^{Q'} \bigl[\bigl(A_t^n-M_{\varepsilon}
\bigr)^+ \bigr]
\\
&\leq& 3\varepsilon.
\end{eqnarray*}
This proves the uniform
integrability of $\{A_t^n\}_{n=1}^{\infty}$
under $\Q'$.
\end{pf}

\subsection{Sequence of discrete-time
martingales with zero limit}\label{S5.3}
For our final approximation result, we construct a sequence of
continuous-time, finite-variation processes that are
martingales when sampled at certain discrete times.
We provide conditions under which this sequence
must converge to zero.

%
%
\begin{definition}\label{D5.6}
A \textit{random partition} $\Pi$ is a set of random times $0=T_0\leq
T_1\leq\cdots\leq T_n$. We set\vspace*{1pt}
$|\Pi|(\omega)\define\sup_{1\leq i\leq n}|T_i(\omega
)-T_{i-1}(\omega)|$. Let $\{\Pi^m\}_{m=1}^{\infty}$ be a sequence of
random partitions, possibly\vspace*{1pt} defined\vadjust{\goodbreak} on different spaces
$\{\Omega^m\}_{m=1}^{\infty}$, where the random times in the partitions
$\Pi^m$ are denoted $T_0^m\leq T_1^m\leq\cdots\leq T_{N(m)}^m$. We say
that $\{\Pi^m\}_{m=1}^{\infty}$ \textit{converges uniformly to the
identity} if
%
%
\begin{equation}
\label{5.13a} \lim_{m\rightarrow\infty}\sup_{\omega\in\Omega^m} \bigl|\Pi^m(
\omega)\bigr|=0
\quad\mbox{and}\quad \lim_{m\rightarrow\infty}\inf_{\omega\in\Omega^m}
T_{N(m)}^m(\omega)=\infty.
\end{equation}
\end{definition}

%
%
\begin{proposition}\label{P5.7}
Let $(\Omega^m,\sF^m,\P^m)_{m=1}^{\infty}$ be a sequence of probability
spaces. Assume that on each space there is defined an
$\R^d$-valued process $X^m$
and a random partition
$\Pi^m=\{T_0^m,T_1^m,\ldots, T_{N(m)}^m\}$,
and these partitions converge uniformly to the identity.
Assume further that the set of processes and measures
$(X^m,\lambda_{[0,t]}\times\P^m)_{m=1}^{\infty}$
is uniformly
integrable for every $t\geq0$.
For $k=0,1,\ldots,N(m)$, define
$Y_k^m\define\int_0^{T_k^m}X_u^m \,du$ and
$\sF_k^m\define\sigma(Y_j^m,T_j^m | 0\leq j\leq k)$,
and assume that $(Y_k^m,\sF_k^m)_{0\leq k\leq N(m)}$ is
a martingale for each $m$. Then
\[
\lim_{m\rightarrow\infty}\E^m \sup_{0\leq s\leq t}\biggl\llVert\int
_0^sX_u^m \,du\biggr
\rrVert=0,\qquad t\geq0.
\]
\end{proposition}

\begin{pf}
By considering components of $\int_0^sX^m_u\,du$, we may reduce
the proof to the case $d=1$. Fix $t\geq0$.
Fix $m$ large enough that ${\sup_{\omega\in\Omega^m}}|\Pi^m|(\omega
)\leq1$
and $\inf_{\omega\in\Omega^m}T_{N(m)}^m(\omega)>t$.
Define $\rho\define\min\{k\dvtx T_k^m\geq t\}$, so that
$T_{\rho}^m$ is the first random time
after $t$ and
$T_{\rho}^m\leq T_{N(m)}^m\wedge(t+1)$. The discrete-time
martingale $Y^m$ stopped at $T_{\rho}^m$ is still a martingale.
For $0\leq s\leq t$, set
$\tau(s)\define\max\{k\dvtx T_k^m\leq s\}$,
so that $T_{\tau(s)}^m$ is the last random time before $s$.
Then $\tau(s)\leq\rho$ and $0\leq s-T_{\tau(s)}^m\leq|\Pi^m|$.
For $M>0$,
\begin{eqnarray*}
\biggl\llvert\int_0^sX_u^m
\,du\biggr\rrvert&\leq& \bigl|Y_{\tau(s)}^m\bigr| +\int
^s_{T_{\tau(s)}^m}\bigl|X^m_u\bigr| \,du
\\
&\leq& \bigl|Y_{\tau(s)}^m\bigr|+\int_{T_{\tau(s)}^m}^s
\bigl[ \bigl(\bigl|X_u^m\bigr|-M \bigr)^++M \bigr]\,du
\\
&\leq& \max_{1\leq k\leq\rho}\bigl|Y_k^m\bigr| +\int
_0^t \bigl(\bigl|X_u^m\bigr|-M
\bigr)^+\,du+M\bigl|\Pi^m\bigr|.
\end{eqnarray*}
Maximizing over $s\in[0,t]$ and taking expectations, we obtain
%
%
\begin{eqnarray}
\label{5.15} \E^m\sup_{s\in[0,t]}\biggl\llvert\int
_0^sX_u^m \,du\biggr
\rrvert&\leq&\E^m\max_{1\leq k\leq\rho}\bigl|Y_k^m\bigr|\nonumber\\
&&{}+\E^m\int_0^t
\bigl(\bigl|X_u^m\bigr|-M \bigr)^+\,du\\
&&{} +M\E^m\bigl|
\Pi^m\bigr|.\nonumber
\end{eqnarray}

We bound the first term on the right-hand side of (\ref{5.15}).
The discrete-time Burkholder--Davis--Gundy
inequality (e.g.,~\cite{Gar}, inequality II.1.1) implies the existence of a
universal constant $C$ for which
%
%
\begin{equation}
\label{5.16} \E^m\max_{1\leq k\leq\rho}\bigl|Y_k^m\bigr|
\leq C\E^m \biggl[ \biggl(\sum_{1\leq k\leq\rho}
\bigl(Y_k^m-Y_{k-1}^m
\bigr)^2 \biggr)^{1/2} \biggr].
\end{equation}
The right-hand side of (\ref{5.16}) can be bounded
using H\"older's inequality. In particular,
%
%
\begin{eqnarray}\label{5.17}
&&
\E^m \biggl[ \biggl(\sum_{1\leq k\leq\rho}
\bigl(Y_k^m-Y_{k-1}^m
\bigr)^2 \biggr)^{1/2} \biggr]
\nonumber
\\
&&\qquad\leq\E^m \biggl[\max_{1\leq k\leq\rho}\bigl|Y_k^m-Y_{k-1}^m\bigr|^{1/2}
\cdot\biggl(\sum_{1\leq k\leq\rho}\bigl|Y_k^m-Y_{k-1}^m\bigr|
\biggr)^{1/2} \biggr]
\nonumber
\\
&&\qquad\leq\E^m \biggl[\max_{1\leq k\leq\rho}\bigl|Y_k^m-Y_{k-1}^m\bigr|^{1/2}
\cdot\biggl(\int_0^{t+1}\bigl|X_u^m\bigr|
\,du \biggr)^{1/2} \biggr]
\\
&&\qquad\leq\sqrt{\E^m\max_{1\leq k\leq\rho}\bigl|Y_k^m-Y_{k-1}^m\bigr|}
\cdot\sqrt{\E^m\int_0^{t+1}\bigl|X_u^m\bigr|
\,du}
\nonumber
\\
&&\qquad\leq\sqrt{\E^m\int_0^{t+1}
\bigl(\bigl|X_u^m\bigr|-M \bigr)^+ \,du+M\E^m\bigl|
\Pi^m\bigr|} \cdot\sqrt{\E^m\int
_0^{t+1}\bigl|X_u^m\bigr| \,du}.\nonumber
\end{eqnarray}
Combining (\ref{5.15})--(\ref{5.17}), we obtain
\begin{eqnarray*}
&&
\E^m\sup_{s\in[0,t]} \biggl\llvert\int
_0^sX^m_u \,du\biggr
\rrvert
\\
&&\qquad\leq C\sqrt{\E^m\int_0^{t+1}
\bigl(\bigl|X_u^m\bigr|-M \bigr)^+ \,du+M\E^m\bigl|
\Pi^m\bigr|} \cdot\sqrt{\E^m\int
_0^{t+1}\bigl|X_u^m\bigr| \,du}
\\
&&\qquad\quad{} +\E^m\int_0^t
\bigl(\bigl|X_u^m\bigr|-M \bigr)^+\,du +M\E^m|\Pi|,
\end{eqnarray*}
where $C$ does not depend on $X$ and $M\geq0$ is arbitrary. The uniform
integrability of $(X^m,\lambda_{[0,t+1]}\times\P^m)_{m=1}^{\infty}$
implies that $\sup_{m}\E^m\int_0^{t+1}|X_u^m|\,du$ is a finite
constant~$C'$. Given $\varepsilon>0$, uniform integrability further
permits us to choose $M$ so large that
$\sup_m\E^m\int_0^{t+1}(|X_u^m|-M)^+\,du\leq\varepsilon$. For such an
$M$,
\[
\E^m\sup_{s\in[0,t]}\biggl\llvert\int_0^sX^m_u
\,du\biggr\rrvert\leq C\sqrt{C'}\sqrt{\varepsilon+ME^m\bigl|
\Pi^m\bigr|}+\varepsilon+M\E^m\bigl|\Pi^m\bigr|.
\]
Letting $m\rightarrow\infty$ and using the first part of (\ref{5.13a}),
we conclude that
\[
\limsup_{m\rightarrow\infty} \E^m\sup_{s\in[0,t]}
\biggl\llvert\int_0^sX^m_u
\,du\biggr\rrvert\leq C\sqrt{C'\varepsilon}+\varepsilon.
\]
\upqed
\end{pf}

\section{\texorpdfstring{Proof of Theorem \protect\ref{T2.4}}{Proof of Theorem 3.6}}
\label{S6}

We prove a theorem that is little
more than a restatement of Theorem~\ref{T2.4} without
reference to the driving Brownian motions $W$ and $\What$
in that theorem. We develop this connection immediately
after the statement of Theorem~\ref{T6.1} below.
Recall Definition~\ref{D3.11}.

%
%
\begin{theorem}\label{T6.1} Let $\sE$ be a Polish space. Let
$(\Omega,\sF,\{\sF_t\}_{t\geq0},\P)$ be a
filtered probability space that supports an
$\sE$-valued random variable $Z_0$ and
an adapted continuous $\R^d$-valued semimartingale
$Y$ with $Y_0=0$ and with characteristic pair
$(B,C)$, where
%
%
\begin{equation}
\label{6.0} B_t=\int_0^tb_s
\,ds,\qquad C_t=\int_0^tc_s
\,ds,
\end{equation}
and the adapted $\R^d$-valued
process $b$ and the adapted $\R^d\times\R^d$-valued
positive semidefinite process $c$
satisfy
%
%
\begin{equation}
\label{6.1} \E\biggl[\int_0^t \bigl(
\|b_s\|+\|c_s\| \bigr) \,ds \biggr]<\infty,\qquad t\geq0.
\end{equation}
Let $\widehat{b}$ and $\widehat{c}$ be measurable
functions defined on $[0,\infty)\times\sE$ with
$\widehat{b}$ taking values in $\R^d$ and $\widehat{c}$ taking
values in the space of $d\times d$ positive
semidefinite matrices, and let $N\subset[0,\infty)$ be a
Lebesgue-null set such that
%
%
\begin{equation}
\label{6.2}\qquad \widehat{b}(t,Z_t)= \E[b_t |
Z_t],\qquad \widehat{c}(t,Z_t)= \E[c_t |
Z_t],\qquad \P\mbox{-a.s.}, t\in N^c.
\end{equation}
Define
$\Omega^{\calE,d}\define\sE\times C_0^d$, let
$\Phi\dvtx\Omega^{\calE,d}\rightarrow C^{\sE}$ be a continuous
updating function and let $Z$ be the continuous,
$\sE$-valued process given by $Z=\Phi(Z_0,Y)$.
Let $\Yhat\dvtx\Omega^{\calE,d}\rightarrow C_0^d$ be
given by $\Yhat(e,x)=x$ and
$\Zhat\dvtx\Omega^{\calE,d}\rightarrow C^{\sE}$ be given
by $\Zhat= \Phi(e,x)$. Then there exists a
measure $\Phat$ on $\Omega^{\calE,d}$ such that:
\begin{longlist}[(ii)]
\item[(i)]
$\Yhat$ is\vspace*{1pt} a semimartingale with
characteristic pair $(\Bhat,\Chat)$ under $\Phat$,
where $\Bhat_t\define\int_0^t\widehat{b}(s,\Zhat_s)\,ds$ and
$\Chat_t\define\int_0^t\widehat{c}(s,\Zhat_s)\,ds$, and
\item[(ii)]
for each $t\geq0$, the distribution of
$\Zhat_t$ under $\Phat$
agrees with the distribution of $Z_t$ under $\P$.
\end{longlist}
\end{theorem}

\begin{pf*}{Proof of Theorem~\ref{T2.4}}
Let us assume Theorem~\ref{T6.1}. Then, under the hypotheses of Theorem
\ref{T2.4}, we may define $c_s\define\sigma_s\sigma_s^{tr}$ and invoke
Proposition~\ref{P4.1} to ensure the existence of functions
$\widehat{b}$ and $\widehat{c}$ and a Lebesgue-null set $N$ such that
(\ref{6.2}) holds. We then conclude that there exist $\Yhat$ and
$\Zhat$ satisfying properties (i) and (ii) in Theorem~\ref{T6.1}. To show
that $\Yhat$ has the representation (\ref{2.14}),\vadjust{\goodbreak} we set $\sigmahat$
equal to the symmetric square root of $\widehat{c}$ and invoke the
It\^o integral representation (e.g.,~\cite{KS}, Chapter 3, Theorem 4.2)
for the $d$-dimensional local martingale $\Yhat-\Bhat$.
\end{pf*}

\begin{pf*}{Proof of Theorem~\ref{T6.1}} %
The proof, which involves a discretization, as suggested by
the example in Section~\ref{DiscreteEx}, and then passage to the
limit, proceeds in several steps.\vspace*{8pt}

\textit{Step} 1: \textit{Construction of canonical space and
processes}.\quad
The random object of interest, $(Z_0,Y,B,C)$, takes values
in $\Omega^{\calE,d}\times C_0^d\times C_0^{d^2}$. In order to show
that the
discretization has a limit, we need to randomize the discretization
times, and thus introduce an extra dimension, defining
$\Omega^*\define[0,1]\times\Omega^{\calE,d}\times C_0^d\times C_0^{d^2}$.
Note that $\Omega^*$ can also be written as
$\Omega^{\sE^*,d^*}$, where $\sE^*=[0,1]\times\sE$
and $d^*=d+d+d^2$. We denote $\sF^{\sE^*,d^*}$ simply
as $\sF^*$ and denote $\sF^{\sE^*,d^*}_t$ simply as $\sF^*_t$.
On $\sF^*$ we define the measure
$\Q$ to be the product of uniform measure on $[0,1]$
and the measure induced by $(Z_0,Y,B,C)$ under $\P$
on $\Omega^{\calE,d}\times C_0^d\times C_0^{d^2}$.
The generic element of $\Omega^*$ will be denoted
$\omega=(\mu,\varepsilon,\eta,\beta,\gamma)$,
and we define the projections
\[
U^*(\omega)=\mu,\qquad Z_0^*(\omega)=\varepsilon,\qquad Y^*(\omega)=\eta,\qquad
B^*(\omega)=\beta,\qquad C^*(\omega)=\gamma.
\]
On the filtered probability space
$(\Omega^*,\sF^*,\{\sF^*_t\}_{t\geq0},\Q)$,
$Y^*$ is a semimartingale with characteristic pair
$(B^*, C^*)$.

We choose an $\R^d$-valued predictable process $b^*$
whose $i$th component at each time $t>0$, denoted
$(b_i^*)_t$, agrees with
\[
\liminf_{k\rightarrow\infty} k \bigl(\bigl(B_i^*\bigr)_t-
\bigl(B^*_i\bigr)_{(t-{1}/{k})^+} \bigr),
\]
whenever the latter is finite.
Likewise, we choose an $\R^{d^2}$-valued predictable process
$c^*$ whose $(i,j)$th component at each time $t>0$,
denoted $(c_{i,j}^*)_t$, agrees with
\[
\liminf_{k\rightarrow\infty} k \bigl(\bigl(C_{i,j}^*\bigr)_t-
\bigl(C^*_{i,j}\bigr)_{(t-{1}/{k})^+} \bigr),
\]
whenever the latter is finite.
By assumption, the components of $B^*$
and $C^*$ are $\Q$-almost surely
absolutely continuous, and so their left derivatives
are defined for Lebesgue-almost every $t\geq0$, $\Q$-almost
surely. By construction,
$b^*$ and $c^*$ are these left derivatives
whenever they are defined. It follows that
%
%
\begin{equation}
\label{6.4}\qquad
\Q\biggl[\int_0^t \bigl(
\bigl\|b^*_s\bigr\|+\bigl\|c^*_s\bigr\| \bigr) \,ds<\infty, B^*_t=
\int_0^tb^*_s \,ds,
C^*_t=\int_0^tc^*_s \,ds
\ \forall t \biggr]=1.
\end{equation}
For $i,j=1,\ldots,d$, the sets
%
%
\begin{eqnarray}\label{6.4a}
&&
\biggl\{t\in[0,\infty)\dvtx\frac{\partial}{\partial t}\bigl(B_i^*
\bigr)_t(\omega) \mbox{ exists but is not equal to }
\bigl(b_i^*\bigr)_t(\omega)\biggr\},
\nonumber\\[-8pt]\\[-8pt]
&&
\biggl\{t\in[0,\infty)\dvtx\frac{\partial}{\partial t}\bigl(C_{i,j}^*
\bigr)_t(\omega) \mbox{ exists but is not equal to }
\bigl(c_{i,j}^*\bigr)_t(\omega)\biggr\}
\nonumber
\end{eqnarray}
are empty for \textit{every} $\omega\in\Omega^*$.\vadjust{\goodbreak}

We set $Z^*\define\Phi(Z_0^*,Y^*)$ and
observe that the random time $U^*$ is strongly independent of
$(Y^*,Z^*,B^*,C^*,b^*,c^*)$ (recall Definition~\ref{D4.3}).
Furthermore, the distribution of $(Y^*,Z^*,B^*,C^*)$ under $\Q$
is the same as the distribution of $(Y,Z,B,C)$ under $\P$, so
(\ref{6.0}) and (\ref{6.4}) imply that
%
%
\begin{equation}
\label{6.45} \E{\int_0^t f(Y_s,
Z_s, b_s, c_s) \,ds} = \E^{\Q} {
\int_0^t f\bigl(Y^*_s,
Z^*_s, b^*_s, c^*_s\bigr) \,ds}
\end{equation}
for any $t \geq0$ and $f$ such that one side of (\ref{6.45}) is well defined.
In particular, (\ref{6.1}) and (\ref{6.45}) imply that
%
%
\begin{equation}
\label{6.5} \E^{\Q} \biggl[\int_0^t
\bigl(\bigl\|b_s^*\bigr\|+\bigl\|c_s^*\bigr\| \bigr) \,ds \biggr]<\infty,
\qquad t\geq0,
\end{equation}
and (\ref{6.2}), (\ref{6.45}) and
Lemma~\ref{L4.2} ensure the existence of a
Lebesgue-null set $N^* \subset[0, \infty)$ such
that
$\widehat{b}(t,Z_t^*)=\E^{\Q}[b_t^* | Z^*_t]$
and $\widehat{c}(t,Z_t^*)=\E^{\Q}[c_t^* | Z^*_t]$
for all $t\notin N^*$.
From (\ref{6.5}) and the conditional version of Jensen's inequality,
we also have
\[
\E^{\Q} \biggl[\int_0^t \bigl(\bigl\|
\widehat{b}\bigl(s,Z_s^*\bigr)\bigr\| +\bigl\|\widehat{c}\bigl(s,Z_s^*
\bigr)\bigr\| \bigr) \,ds \biggr]<\infty,\qquad t\geq0,
\]
or equivalently,
\[
\E\biggl[\int_0^t \bigl(\bigl\|
\widehat{b}(s,Z_s)\bigr\| +\bigl\|\widehat{c}(s,Z_s)\bigr\| \bigr) \,ds
\biggr]<\infty,\qquad t\geq0.
\]

\textit{Step} 2: \textit{Construction of extended partitions}.\quad
For each positive integer $m$, set $N(m)=m^2$,
$T_0^m\define0$ and for $i=1,\ldots,N(m)$, set
$T_i^m\define(U^*+i-1)/m$. Note that each $T_i^m$
is $\sigma(U^*)$ measurable, and consequently
is an $\{\sF^*_t\}_{t\geq0}$-stopping time. Let $\Pi^m$
denote this set of stopping times. The sequence
of random partitions $\{\Pi^m\}_{m=1}^{\infty}$
converges uniformly to the identity (Definition~\ref{D5.6}).

For the next step, we adopt the notation $X=(Y^*,B^*,C^*)$. We set
$\sG_0^m=\sH_0^m=\sF_0^*=\sigma(U,Z_0^*)$, and for $i=1,\ldots,N(m)$,
we set $\sG_i^m=\sigma(U^*,Z^*_{T_i})$ and
$\sH_i^m=\sG_{i-1}^m\vee\sigma( \Delta(X^{T_i^m},T_{i-1}^m) )$.
Finally,\vspace*{2pt} we set $T_{N(m)+1}^m=\infty$ and $\sH_{N(m)+1}^m=
\sG_{N(m)}^m\vee\sigma (\Delta(X,T_{N(m)}^m) )$. It is clear that part
(a) of Definition~\ref{D3.2} is satisfied. To show that
$(T_i^m,\sG_i^m)_{i=1}^{N(m)}$ is an extended partition, it suffices to
verify condition (b) of Definition~\ref{D3.2}, that is, that
$\sG_i^m\subset\sH_i^m$ for $i=1,\ldots,N(m)$. In particular, it
suffices to show that $Z^*_{T_i^m}$ is measurable with\vspace*{-1pt} respect to
$\sigma(U^*)\vee\sigma(Z^*_{T_{i-1}^m})
\vee\sigma(\Delta(X^{T_i^m},T_{i-1}^m))$. Let $\tau\geq0$ be a possibly
random time and define $S_i^m=(\tau-T_{i-1}^m)^+$. On the set
$T_{i-1}^m\leq\tau\leq T_i^m$, we may use property (\ref{2.7}) of the
updating function $\Phi$ to write
%
%
\begin{eqnarray}\label{6.10a}
Z^*_{\tau} &=& \Theta_{S_i^m} \bigl(Z^*,T_{i-1}^m
\bigr)
\nonumber
\\
&=& \Theta_{S_i^m} \bigl(\Phi\bigl(Z_0^*,Y^*
\bigr),T_{i-1}^m \bigr)
\nonumber\\[-8pt]\\[-8pt]
&=& \Phi_{S_i^m} \bigl(\Phi_{T_{i-1}^m}\bigl(Z_0^*,Y^*
\bigr), \Delta\bigl(Y^*,T_{i-1}^m\bigr) \bigr)
\nonumber
\\
&=& \Phi_{S_i^m} \bigl(Z^*_{T_{i-1}^m}, \Delta\bigl(Y^*,T_{i-1}^m
\bigr) \bigr).\nonumber
\end{eqnarray}
If we take $\tau=T_i^m$, this leads to
\[
Z_{T_i^m}^*= \Phi_{T_i^m-T_{i-1}^m} \bigl(Z^*_{T_{i-1}^m}, \Delta
\bigl(Y^*,T_{i-1}^m\bigr) \bigr) \\
=\Phi_{T_i^m-T_{i-1}^m}^{T_i^m-T_{i-1}^m}
\bigl(Z^*_{T_{i-1}^m}, \Delta\bigl(Y^*,T_{i-1}^m\bigr)
\bigr),
\]
and by property (\ref{2.6}), the last expression depends on
the path of $\Delta(Y^*,T_{i-1}^m)$ only up to
time $T_i^m-T_{i-1}^m$, which agrees with the path
of $\Delta((Y^*)^{T_i^m},T_{i-1}^m)$ up to time
$T_i^m-T_{i-1}^m$.
We have thus written $Z^*_{T_i^m}$ in terms of
$T_i^m-T_{i-1}^m$, which is nonrandom unless $i=1$, in which
case it is $U^*/m$, in terms of $Z^*_{T_{i-1}^m}$,
and in terms of $\Delta((Y^*)^{T_i^m},T_{i-1}^m)$.

\textit{Step} 3: \textit{Concatenated measures}.\quad We denote by
$\Pi^m$ the extended partition $(T_i^m,\sG_i^m)_{i=1}^{N(m)}$. These
extended partitions are on the space
$\Omega^*=[0,1]\times\Omega^{\calE,d}\times C_0^d\times C_0^{d^2}$,
which is the same as $\Omega^{\sE^*,d^*}$ defined in step 1. Theorem
\ref{T3.3} implies the existence of concatenated measures
$\Q^m\define\Q^{\otimes\Pi^m}$ that satisfy
%
%
\begin{eqnarray}
\label{6.11b}
\Q^m[A]&=&\Q[A],\qquad A\in\sH_i^m,\qquad i=0,1,
\ldots,N(m)+1,
\\
\label{6.11c}
\qquad\Q^m\bigl[B | \sF^*_{T_i^m}\bigr] &=& \Q\bigl[B |
\sG_i^m\bigr],\qquad B\in\sH_{i+1}^m,\qquad
i=0,1,\ldots,N(m).
\end{eqnarray}
Applying Proposition~\ref{P3.12} with $X=(Y^*,B^*,C^*)$,
we see that $Y^*$ is a semimartingale with
characteristic pair $(B^*,C^*)$ under each $\Q^m$.\vspace*{8pt}

\textit{Step} 4: \textit{Tightness and convergence}.\quad Corollary
\ref{C3.11}(b) shows that the collection of measures induced on
$C_0^d\times C_0^{d^2}$ by $(B^*,C^*)$ under $\{\Q^m\}_{m=1}^{\infty}$
is tight. Theorem~VI.4.18 of~\cite{JS} (Rebolledo's criterion; see
\cite{R}) then implies that the collection of measures induced on
$C_0^d$ by $Y^*$ under $\{\Q^m\}_{m=1}^{\infty}$ is tight. Since
$Z_0^*$ has the same distribution under every $\Q^m$, the set of
measures induced on $\Omega^{\calE,d}$ by $(Z_0^*,Y^*)$ is likewise
tight. Passing to a convergent subsequence if necessary, we obtain a
limiting measure $\Phat$ on $\Omega^{\calE,d}$. To simplify notation,
we assume that the passage to a subsequence is not necessary to
obtain\vspace*{1pt} convergence. We denote the coordinate mappings on
$\Omega^{\calE,d}$ by $\Zhat_0$ and $\Yhat$, and we define
$\Zhat=\Phi(\Zhat_0,\Yhat)$. The continuous mapping theorem implies
that the distributions of $(Y^*,Z^*)$ on $C_0^d\times C^{\sE}$ under
the sequence of measures $\{\Q^m\}_{m=1}^{\infty}$
converge\vspace*{2pt} to the distribution of $(\Yhat,\Zhat)$ under
$\Phat$, that is, $\Q^m\circ(Y^*,Z^*)^{-1}
\Longrightarrow\Phat\circ(\Yhat,\Zhat)^{-1}$.\vspace*{8pt}

\textit{Step} 5: \textit{Agreement of one-dimensional
distributions}.\quad
Returning to (\ref{6.10a}), we take $\tau=t$, a fixed nonnegative
number, so that $S_i^m=(t-T_{i-1}^m)^+$.
On the $\sH_i^m$-measurable
set $\{T_{i-1}^m\leq t<T_i^m\}$, we have
\[
Z_t^*= \Phi_{S_i^m} \bigl(Z^*_{T_{i-1}^m},\Delta
\bigl(Y^*,T_{i-1}^m\bigr) \bigr),
\]
and the term
$\Phi_{S_i^m} (Z^*_{T_{i-1}^m},\Delta(Y^*,T_{i-1}^m) )$
restricted to $\{T_{i-1}^m\leq t<T_i^m\}$
depends only on $S_i^m$, $Z^*_{T_{i-1}^m}$ and
$\Delta((Y^*)^{T_i^m},T_{i-1}^m)$, all of which are
$\sH_i^m$-measurable. Because $\Q^m$ and $\Q$ agree
on each $\sH_i^m$, we conclude that for every Borel
subset $A$ of $\sE$ and for every $t\geq0$,
%
%
\begin{eqnarray}\label{6.15a}
\Q^m\bigl[Z_t^*\in A\bigr] &=& \sum
_{i=1}^{N(m)+1}\Q^m\bigl[Z^*_t
\in A\mbox{ and }T_{i-1}^m\leq t<T_i^m
\bigr]
\nonumber
\\
&=& \sum_{i=1}^{N(m)+1} \Q
\bigl[Z^*_t\in A\mbox{ and }T_{i-1}^m\leq
t<T_i^m\bigr]
\\
&=& \Q\bigl[Z_t^*\in A\bigr] = \P[Z_t\in A].\nonumber
\end{eqnarray}
But the distributions of $Z^*$ under the sequence of measures
$\{\Q^m\}_{m=1}^{\infty}$ converge to the distribution
of $\Zhat$ under $\Phat$, and part (ii) of Theorem~\ref{T6.1}
is proved.\vspace*{8pt}

\textit{Step} 6: \textit{Semimartingale characteristics of the
limit}.\quad
To complete the proof, we must show that under the
measure $\Phat$ on $\Omega^{\calE,d}$,
$\Yhat$ is a semimartingale with characteristic pair
$(\Bhat,\Chat)$, defined in part (i) of Theorem~\ref{T6.1}.
We do this by showing that the distribution
of the $(Y^*,Z^*,B^*,C^*)$ under $\Q^m$ converges to
the distribution of $(\Yhat,\Zhat,\Bhat,\Chat)$ under $\Phat$,
that is,
%
%
\begin{equation}
\label{6.11} \Q^m\circ\bigl(Y^*,Z^*,B^*,C^*\bigr)^{-1}
\Longrightarrow\Phat\circ(\Yhat,\Zhat,\Bhat,\Chat)^{-1}.
\end{equation}
The filtration on $\Omega^{\calE,d}$, defined at the beginning of
Section~\ref{SConMea}, is generated by~$\Yhat$.
Once (\ref{6.11}) is established,
Theorem IX.2.4 of~\cite{JS} will give the desired result.

On $\Omega^*$ we define the processes
\begin{eqnarray*}
\overline{b}_t&\define&\widehat{b}\bigl(t,Z_t^*\bigr),\qquad
\Bbar_t\define\int_0^t
\overline{b}_s \,ds, \\
\overline{c}_t&\define&\widehat{c}
\bigl(t,Z_t^*\bigr), \qquad\Cbar_t\define\int
_0^t\overline{c}_s \,ds,\qquad t\geq0.
\end{eqnarray*}
According to Proposition~\ref{P5.2},
%
%
\begin{eqnarray}
\label{6.12}\quad
&\displaystyle \Q^m\circ\bigl(Y^*,Z^*,\Bbar,\Cbar\bigr)^{-1}
\Longrightarrow\Phat\circ(\Yhat,\Zhat,\Bhat,\Chat)^{-1},&
\\
\label{6.13}
&\displaystyle \bigl\{\overline{b},\lambda_{[0,t]}\times\Q^m\bigr
\}_{m\in\N}\mbox{ is uniformly integrable for every }t\in[0,\infty),&
\nonumber\\[-8pt]\\[-8pt]
&\displaystyle \bigl\{\overline{c},\lambda_{[0,t]}\times\Q^m\bigr
\}_{m\in\N}\mbox{ is uniformly integrable for every }t\in[0,\infty).&
\nonumber
\end{eqnarray}
We show that
$\Q^m\circ(Y^*,Z^*,B^*,C^*)^{-1}$
and $\Q^m\circ(Y^*,Z^*,\Bbar,\Cbar)^{-1}$
have the same limit as $m\rightarrow\infty$.
We do this by showing that
for every $\varepsilon>0$ and $t\in[0,\infty)$,
%
%
\begin{eqnarray}
\label{6.16}
\lim_{m\rightarrow\infty}\Q^m \Bigl[\sup_{0\leq s\leq t}
\bigl\|B^*_s-\Bbar_s\bigr\|\geq\varepsilon\Bigr]&=&0,
\\
\label{6.17}
\lim_{m\rightarrow\infty}\Q^m \Bigl[\sup_{0\leq s\leq t}
\bigl\|C^*_s-\Cbar_s\bigr\|\geq\varepsilon\Bigr]&=&0.
\end{eqnarray}
Once this has been done, (\ref{6.12}) will imply (\ref{6.11}),
and we will be done.\vspace*{8pt}

\textit{Step} 7: \textit{Proof of} (\ref{6.16}) \textit{and}
(\ref{6.17}).\quad
In fact, we prove only (\ref{6.16}), because the
proof of (\ref{6.17}) is the same. Without loss of
generality, we assume that $B^*$ and $\Bhat$ are one dimensional.

For $i=1,\ldots,N(m)$, define the
$\sH_{i+1}^m$-measurable random variable
\[
\xi_i^m\define\liminf_{k\rightarrow\infty} k
\bigl(B^*_{T_i^m+{1}/{k}}-B^*_{T_i^m}\bigr),
\]
which is the right derivative of $B^*$ at $T_i^m$ whenever
this derivative is defined. Recall from step 1
that $b^*_{T_i^m}$ is the left
derivative of $B^*$ at $T_i^m$ whenever this derivative
is defined and is finite. By construction, $B^*$ is independent of
$T_i^m$ under $\Q$, and its derivative is defined and is finite
Lebesgue-almost everywhere, $\Q$-almost surely.
But $T_i^m$ is uniformly distributed on
$[\frac{i-1}{n},\frac{i}{n}]$. It follows that
%
%
\begin{equation}
\label{6.19} \Q\bigl[\xi_i^m=b^*_{T_i^m}
\bigr]=1,\qquad i=1,\ldots,N(m).
\end{equation}

We define three sequences of step functions:
\begin{eqnarray*}
b_t^m&\define&\sum_{i=1}^{N(m)}
\xi_i^m\ind_{[T_i^m,T_{i+1}^m)}(t),\qquad\overline{b}{}^m_t
\define\sum_{i=1}^{N(m)}\overline{b}_{T_i^m}
\ind_{[T_i^m,T_{i+1}^m)},\\
b_t^{\Pi^m}&\define&\sum
_{i=1}^{N(m)}b_{T_i^n}^* \ind_{[T_i^m,T_{i+1}^m)}.
\end{eqnarray*}
We further define
\[
B_t^m\define\int_0^tb_s^m
\,ds,\qquad \Bbar{}^m_t\define\int_0^t
\overline{b}{}^m_s \,ds.
\]
Because of (\ref{6.19}), $b^m$ and $b^{\Pi^m}$ are
$\Q$-indistinguishable.

Each $B^m$ is piecewise linear,
and so for \textit{every} $\omega\in\Omega^*$,
$\frac{\partial}{\partial t}B_t^m(\omega)=b_t^m(\omega)$
except at finitely many values of $t$.
In addition, $\Delta(B^m,T_i^m)$ is
$\sigma(\xi_j^m\dvtx j\geq i)$-measurable.
For $j\geq i$, $\xi_j^m$ is
$\sH_{j+1}^m$-measurable,
and we have shown in the proof of Theorem
\ref{T3.3} [see (\ref{HinG})]
that $\sH_{j+1}^m\subset\sG_i^m\vee\sigma(\Delta
(X,T_i^m) )$
for $j=i,i+1,\ldots,N(m)$, so we may conclude that
$\Delta(B^m,T_i^m)$ is
$\sG_i^m\vee\sigma(\Delta(X,T_i^m) )$-measurable
for $i=1,\ldots,N(m)$.
This measurability condition
is trivially satisfied when $i=0$ as well.
We conclude that the pair of processes
$(B^m, b^m)$ satisfies the hypotheses
of Proposition~\ref{P3.10}, including (\ref{3.15}),
with $\Q$ replacing $\P$ and $\Q^m$ replacing $\P^{\otimes\Pi}$.

Because $B^*$ is a component of $X$ and the set
(\ref{6.4a}) is empty, $(B^*, b^*)$ also satisfies
the hypothesis of Proposition~\ref{P3.10}, and hence,
so does $(B^*-B^m, b^*-b^m)$. We thus obtain from (\ref{3.17})
that
\[
\E^{\Q^m}\int_0^t\bigl|b_s^*-b_s^m\bigr|
\,ds =\E^{\Q}\int_0^t\bigl|b_s^*-b_s^m\bigr|
\,ds,\qquad t\geq0.
\]
For fixed $t\geq0$, we use this equality,
the $\Q$-indistinguishability of $b^m$ and $b^{\Pi^m}$
and Proposition~\ref{P5.5} to write
%
%
\begin{eqnarray}\label{6.21}
\limsup_{m\rightarrow\infty}\E^{\Q^m} \sup_{0\leq s\leq t}\bigl|B_s^*-B_s^m\bigr|
&\leq& \limsup_{m\rightarrow\infty}\E^{\Q^m} \int_0^t\bigl|b_s^*-b_s^m\bigr|
\,ds
\nonumber
\\
&=& \limsup_{m\rightarrow\infty} \E^{\Q}\int_0^t\bigl|b_s^*-b_s^m\bigr|
\,ds
\nonumber\\[-8pt]\\[-8pt]
&=& \limsup_{m\rightarrow\infty} \E^{\Q}\int_0^t\bigl|b_s^*-b_s^{\Pi^m}\bigr|
\,ds
\nonumber
\\
&=& 0.\nonumber
\end{eqnarray}

We consider the difference between $\Bbar$
and $\Bbar{}^m$. For $i=1,\ldots,N(m)+1$,
%
%
\begin{eqnarray} \label{6.22}
\quad&&
\int_{T_{i-1}^m\wedge t}^{T_i^m\wedge t} \bigl|\overline{b}_s-
\overline{b}{}^m_s\bigr| \,ds
\nonumber
\\
&&\qquad= \int_{T_{i-1}^m\wedge t}^{T_i^m\wedge t} \bigl|\widehat{b}
\bigl(s,Z_s^*\bigr)-\widehat{b}\bigl(T_{i-1}^m,Z^*_{T_{i-1}^m}
\bigr)\bigr| \,ds
\\
&&\qquad= \int_0^{S_i^m} \bigl|\widehat{b}
\bigl(T_{i-1}^m+s,\Theta_s \bigl(\bigl(Z^*
\bigr)^{T_i^m},T_{i-1}^m \bigr) \bigr) -\widehat{b}
\bigl(\Theta_0 \bigl(\bigl(Z^*\bigr)^{T_i^m},T_{i-1}^m
\bigr) \bigr) \bigr| \,ds,\nonumber
\end{eqnarray}
where $S_i^m=\frac{1}{n}\wedge(t-T_{i-1}^m)^+$ if $i\geq2$
and $S_1^m=T_1^m\wedge t$.
The final expression in (\ref{6.22}) is $\sH_i^m$-measurable,
and so the first expression is as well.
But $\Q^m$ and $\Q$ agree on $\sH_i^m$, which together
with Proposition~\ref{P5.5} implies
%
%
\begin{eqnarray}\label{6.23}
\limsup_{m\rightarrow\infty}\E^{\Q^m}\sup_{0\leq s\leq t} \bigl|
\Bbar_s-\Bbar{}^m_s\bigr| &\leq&
\limsup_{m\rightarrow\infty} \E^{\Q^m}\int_0^t\bigl|
\overline{b}_s-\overline{b}{}^m_s\bigr| \,ds
\nonumber
\\
&=& \limsup_{m\rightarrow\infty} \sum_{i=1}^{N(m)+1}
\E^{\Q^m} \int_{T_{i-1}^m\wedge t}^{T_i^m\wedge t} \bigl|
\overline{b}_s-\overline{b}{}^m_s\bigr| \,ds
\nonumber\\[-8pt]\\[-8pt]
&=& \limsup_{m\rightarrow\infty} \sum_{i=1}^{N(m)+1}
\E^{\Q} \int_{T_{i-1}^m\wedge t}^{T_i^m\wedge t} \bigl|
\overline{b}_s-\overline{b}{}^m_s\bigr| \,ds
\nonumber
\\
&=& \limsup_{m\rightarrow\infty} \E^{\Q}\int_0^t\bigl|
\overline{b}_s-\overline{b}{}^m_s\bigr| \,ds
\nonumber
\\
&=&0.\nonumber
\end{eqnarray}

It remains to estimate the difference between $B^m$
and $\Bbar{}^m$. From (\ref{6.13}) and (\ref{6.23})
we see that $\{\overline{b}{}^m,\lambda_{[0,t]}\times\Q^m\}_{m\in\N}$
is uniformly integrable for every $t\in[0,\infty)$.
We show that $\{b^m,\lambda_{[0,t]}\times\Q^m\}_{m\in\N}$
is also uniformly integrable by using the $\sH_{i+1}^m$
measurability of $T_i^m$, $T_{i+1}^m$ and $\xi_i^m$ to write
%
%
\begin{eqnarray}\label{6.24}
\E^{\Q^m}\int_0^t\bigl|b_s^m\bigr|
\ind_{\{|b_s^n|\geq M\}} \,ds &=& \sum_{i=0}^{N(m)}
\E^{\Q^m} \bigl[\bigl(T_{i+1}^m\wedge
t-T_i^m\wedge t\bigr) \bigl|\xi_i^m\bigr|
\ind_{\{|\xi_i^m|\geq M\}} \bigr]
\nonumber
\\
&=& \sum_{i=0}^{N(m)}\E^{\Q}
\bigl[\bigl(T_{i+1}^m\wedge t-T_i^m
\wedge t\bigr) \bigl|\xi_i^m\bigr|\ind_{\{|\xi_i^m|\geq M\}} \bigr]
\\
&=& \E^{\Q}\int_0^t\bigl|b_s^m\bigr|
\ind_{\{|b_s^n|\geq M\}} \,ds.\nonumber
\end{eqnarray}
Under $\lambda_{[0,t]}\times\Q$, $b^*$ restricted to $[0,t]$
is integrable [see (\ref{6.4})].
One consequence of (\ref{6.21}) is that $b^m$ restricted
to $[0,t]$ converges
to $b^*$ restricted to $[0,t]$ in $L^1(\lambda_{[0,t]}\times\Q)$.
This, combined with (\ref{6.24}), yields the uniform integrability
of $\{b^m,\lambda_{[0,t]}\times\Q^m\}_{m\in\N}$.
We conclude that
$\{b^m-\overline{b}{}^m,\lambda_{[0,t]}\times\Q^m\}_{m\in\N}$
is uniformly integrable for every $t\geq0$.

Define
\[
\Psi_k^m\define B_{T_k^m}-
\Bbar_{T_k^m} =\int_0^{T_k^m}
\bigl(b_s^m-\overline{b}{}^m_s
\bigr) \,ds.
\]
Let $k=0,1,\ldots,N(m)-1$ be given.
Because $T_{k+1}^m$ and $ T_{k}^m$ are $\sF^*_{T_k^m}$-measurable,
$\xi_{k}^m-\overline{b}_{T_{k}^m}
=\xi_{k}^m-\widehat{b}(T_{k}^m,Z_{T_{k}^m})$ is $\sH_{k+1}^m$-measurable,
and (\ref{6.11c}) and (\ref{6.19}) hold,
we may write
\begin{eqnarray*}
&&
\E^{\Q^m}\bigl[\Psi_{k+1}^m-
\Psi_k^m | \sF^*_{T_{k}^m}\bigr]
\\
&&\qquad= \bigl(T_{k+1}^m-T_{k}^m
\bigr)\E^{\Q^m}\bigl[\xi_{k}^m -\widehat{b}
\bigl({T_{k}^m},Z^*_{T_{k}^m}\bigr) |
\sF^*_{T_{k}^m}\bigr]
\\
&&\qquad= \bigl(T_{k+1}^m-T_{k}^m
\bigr)\E^{\Q}\bigl[\xi_{k}^m -\widehat{b}
\bigl({T_{k}^m},Z^*_{T_{k}^m}\bigr) |
\sG_{k}^m\bigr]
\\
&&\qquad= \bigl(T_{k+1}^m-T_{k}^m
\bigr) \bigl(\E^{\Q}\bigl[b_{T_k^m}^* | \sG_k^m
\bigr] -\widehat{b}\bigl({T_{k}^m},Z^*_{T_{k}^m}
\bigr) \bigr).
\end{eqnarray*}
Proposition~\ref{P4.5} implies that
\[
\E^{\Q}\bigl[b_{T_k^m}^* | \sG_k^m
\bigr] =\E^{\Q}\bigl[b_{T_k^m}^* | T_k^m,Z^*_{T_k^m}
\bigr] =\widehat{b}\bigl(T_k^m,Z^*_{T_k^m}
\bigr).
\]
We conclude that $(\Psi_k^m,\sF^*_{T_k^m} | 0\leq k\leq N(m))$
is a discrete-time martingale under~$\Q^m$, which implies
that $(\Psi_k^m,\sF_k^m | 0\leq k\leq N(m))$ is also
a martingale, where
$\sF_k^m\define\sigma(\Psi_j^m,T_j^m | 0\leq j\leq k)
\subset\sF^*_{T_k^m}$. Proposition~\ref{P5.7} now implies
that
%
%
\begin{equation}
\label{6.29} \lim_{m\rightarrow\infty}\E^{\Q^m}\sup_{0\leq s\leq t}
\bigl|B^m_s-\Bbar{}^m_s\bigr|=0,\qquad t\geq0.
\end{equation}

Using the triangle inequality, we combine (\ref{6.21}),
(\ref{6.29}) and (\ref{6.23}) to conclude
%
%
\begin{equation}
\label{6.30} \limsup_{m\rightarrow\infty}\E^{\Q^m}\sup_{0\leq s\leq t}
\bigl|B^*_s-\Bbar_s\bigr|=0.
\end{equation}
Equation (\ref{6.16}) follows.
\end{pf*}

\section*{Acknowledgments}

We thank Peter Carr for pointing out Gy\"ongy~\cite{G} and an anonymous
referee for a number of helpful comments.


%

\printaddresses

\end{document}